\pdfoutput=1
\RequirePackage{ifpdf}
\ifpdf 
\documentclass[pdftex]{sigma}
\else
\documentclass{sigma}
\fi

\begin{document}

\allowdisplaybreaks

\renewcommand{\PaperNumber}{007}

\FirstPageHeading

\ShortArticleName{$B_{2}$ Matrix Weight Function}

\ArticleName{Vector-Valued Polynomials\\ and a Matrix Weight Function with $\boldsymbol{B_{2}}$-Action}

\Author{Charles F.~DUNKL}

\AuthorNameForHeading{C.F.~Dunkl}

\Address{Department of Mathematics, University of Virginia,\\ PO Box 400137,
Charlottesville VA 22904-4137, USA}
\Email{\href{mailto:cfd5z@virginia.edu}{cfd5z@virginia.edu}}
\URLaddress{\url{http://people.virginia.edu/~cfd5z/home.html}}

\ArticleDates{Received October 16, 2012, in f\/inal form January 23, 2013; Published online January 30, 2013}

\Abstract{The structure of orthogonal polynomials on $\mathbb{R}^{2}$ with the weight
function $\vert x_{1}^{2}-x_{2}^{2}\vert ^{2k_{0}}\vert
x_{1}x_{2}\vert ^{2k_{1}}e^{-(  x_{1}^{2}+x_{2}^{2})  /2}$ is
based on the Dunkl operators of type $B_{2}$. This refers to the full symmetry
group of the square, generated by ref\/lections in the lines $x_{1}=0$ and
$x_{1}-x_{2}=0$. The weight function is integrable if $k_{0},k_{1},k_{0}
+k_{1}>-\frac{1}{2}$. Dunkl operators can be def\/ined for polynomials taking
values in a module of the associated ref\/lection group, that is, a vector space
on which the group has an irreducible representation. The unique
$2$-dimensional representation of the group $B_{2}$ is used here. The specif\/ic
operators for this group and an analysis of the inner products on the harmonic
vector-valued polynomials are presented in this paper. An orthogonal basis for
the harmonic polynomials is constructed, and is used to def\/ine an
exponential-type kernel. In contrast to the ordinary scalar case the inner
product structure is positive only when $(  k_{0},k_{1})  $ satisfy
$-\frac{1}{2}<k_{0}\pm k_{1}<\frac{1}{2}$. For vector polynomials $(
f_{i})  _{i=1}^{2}$, $( g_{i})  _{i=1}^{2}$ the inner product
has the form $\iint_{\mathbb{R}^{2}}f(x)  K(x)
g(x)  ^{T}e^{-(  x_{1}^{2}+x_{2}^{2})  /2}dx_{1}dx_{2}$ where the matrix function $K(x)$ has to satisfy various
transformation and boundary conditions. The matrix~$K$ is expressed in terms
of hypergeometric functions.}

\Keywords{matrix Gaussian weight function; harmonic polynomials}

\Classification{33C52; 42C05; 33C05}

\vspace{-2mm}

\section{Introduction}

The algebra of operators on polynomials generated by multiplication and the
Dunkl operators associated with some ref\/lection group is called the rational
Cherednik algebra. It is parametrized by a multiplicity function which is
def\/ined on the set of roots of the group and is invariant under the group
action. For scalar-valued polynomials there exists a Gaussian-type weight
function which demonstrates the positivity of a certain bilinear form on
polynomials, for positive values (and a small interval of negative values) of
the multiplicity function. The algebra can also be represented on polynomials
with values in an irreducible module of the group. In this case the problem of
f\/inding a Gaussian-type weight function and the multiplicity-function values
for which it is positive and integrable becomes much more complicated. Here we
initiate the study of this problem on the smallest two-parameter,
two-dimensional example, namely, the group of type~$B_{2}$ (the full symmetry
group of the square).

Grif\/feth \cite{Griffeth2010} def\/ined and studied analogues of nonsymmetric
Jack polynomials for arbitrary irreducible representations of the complex
ref\/lection groups in the family~$G(  r,1,n)  $. This paper
introduced many useful methods for dealing with vector-valued polynomials. In
the present paper we consider~$B_{2}$, which is the member $G(
2,1,2)  $ of the family, but we use harmonic polynomials, rather than
Grif\/feth's Jack polynomials because the former play a crucial part in the
analysis of the Gaussian weight. There is a detailed study of the unitary
representations of the rational Cherednik algebra for the symmetric and
dihedral groups in Etingof and Stoica~\cite{Etingof/Stoica:2009}.

We begin with a brief discussion of vector-valued polynomials and the results
which hold for any real ref\/lection group. This includes the def\/inition of the
Dunkl operators and the basic bilinear form. The next section specializes to
the group $B_{2}$ and contains the construction of an orthogonal basis of
harmonic homogeneous polynomials, also a brief discussion of the radical.
Section~\ref{section4} uses this explicit basis to construct the appropriate analogue of
the exponential function. Section~\ref{section5} contains the derivation of the
Gaussian-type weight function; it is a $2\times2$ matrix function whose
entries involve hypergeometric functions. This is much more complicated than
the scalar case. The method of solution is to set up a system of dif\/ferential
equations, f\/ind a fundamental solution and then impose several geometric
conditions, involving behavior on the mirrors (the walls of the fundamental
region of the group) to construct the desired solution.

\section{General results}\label{section2}

Suppose $R$ is a root system in $\mathbb{R}^{N}$ and $W=W(  R)  $
is the f\/inite Coxeter group generated by $\{  \sigma_{v}:v\in
R_{+}\}  $ (where $\langle x,y\rangle :=\sum\limits_{i=1}^{N}
x_{i}y_{i}, \vert x \vert := \langle x,x \rangle
^{1/2}$, $x\sigma_{v}:=x-2\frac{ \langle x,v \rangle }{ \langle
v,v \rangle }v$ for $x,y,v\in\mathbb{R}^{N}$ and $v\neq0$). Let $\kappa$
be a multiplicity function on~$R$ ($u=vw$ for some $w\in W$ and $u,v\in R$
implies $\kappa(u)  =\kappa(v)  $). Suppose $\tau$
is an irreducible representation of~$W$ on a (real) vector space~$V$ of
polynomials in $t\in\mathbb{R}^{N}$ of dimension~$n_{\tau}$. (There is a
general result for these groups that real representations suf\/f\/ice, see
\cite[Chapter~11]{Carter1993}.) Let $\mathcal{P}_{V}$ be the space of polynomial
functions $\mathbb{R}^{N}\rightarrow V$, that is, the generic $f\in
\mathcal{P}_{V}$ can be expressed as $f(x,t)$ where~$f$ is a
polynomial in~$x$,~$t$ and $f(x,t)  \in V$ for each f\/ixed
$x\in\mathbb{R}^{N}$. There is an action of $W$ on $\mathcal{P}_{V}$ given by
\[
wf(x,t)  :=f (  xw,tw )  , \qquad w\in W.
\]
Def\/ine Dunkl operators on $\mathcal{P}_{V}$, for $1\leq i\leq N$ by%
\[
\mathcal{D}_{i}f(x,t)  :=\frac{\partial}{\partial x_{i}}f (x,t)  +\sum_{v\in R_{+}}\kappa(v)  \frac{f(x,t\sigma_{v})  -f(x\sigma_{v},t\sigma_{v})  }{\langle x,v \rangle }v_{i}.
\]
There is an equivariance relation, for $u\in\mathbb{R}^{N}$, $w\in W$%
\begin{gather}
\sum_{i=1}^{N}u_{i}\mathcal{D}_{i}w=w\sum_{i=1}^{N}(uw)
_{i}\mathcal{D}_{i}. \label{wDw}
\end{gather}
Ordinary (scalar) polynomials act by multiplication on $\mathcal{P}_{V}$. For
$1\leq i,j\leq N$ and $f\in\mathcal{P}_{V}$ the basic commutation rule is
\begin{gather}
\mathcal{D}_{i}x_{j}f(x,t)  -x_{j}\mathcal{D}_{i}f (
x,t )  =\delta_{ij}f(x,t)  +2\sum_{v\in R_{+}}\kappa (
v )  \frac{v_{i}v_{j}}{\vert v\vert ^{2}}f (  x\sigma
_{v},t\sigma_{v} )  . \label{dxxd}
\end{gather}

The abstract algebra generated by $ \{  x_{i},\mathcal{D}_{i}:1\leq i\leq
N \}  \cup\mathbb{R}W$ with the commutation relations $x_{i}x_{j}
=x_{j}x_{i}$, $\mathcal{D}_{i}\mathcal{D}_{j}=\mathcal{D}_{j}\mathcal{D}_{i}$,
(\ref{dxxd}) and equivariance relations like~(\ref{wDw}) is called the
\textit{rational Cherednik algebra} of~$W$ parametrized by $\kappa$;
henceforth denoted by $\mathcal{A}_{\kappa}$. Then~$\mathcal{P}_{V}$ is called
the \textit{standard module} of $\mathcal{A}_{\kappa}$ determined by the
$W$-module~$V$.

We introduce symmetric bilinear $W$-invariant forms on $\mathcal{P}_{V}$.
There is a $W$-invariant form $ \langle \cdot,\cdot \rangle _{\tau}$
on~$V$; it is unique up to multiplication by a constant because $\tau$ is
irreducible. The form is extended to $\mathcal{P}_{V}$ subject to
$ \langle x_{i}f(x,t)  ,g(x,t)   \rangle
_{\tau}= \langle f(x,t)  ,\mathcal{D}_{i}g(x,t)
 \rangle _{\tau}$ for $f,g\in\mathcal{P}_{V}$ and $1\leq i\leq N$. To be
more specif\/ic: let $ \{  \xi_{i}(t)  :1\leq i\leq n_{\tau
} \}  $ be a basis for $V$. Any $f\in\mathcal{P}_{V}$ has a~unique
expression $f(x,t)  =\sum_{i}f_{i}(x)  \xi
_{i}(t)  $ where each $f_{i}(x)  $ is a polynomial,
and then{\samepage
\[
\left\langle f,g\right\rangle _{\tau}:=\sum_{i}\left\langle \xi_{i}\left(
t\right)  ,\left(  f_{i}\left(  \mathcal{D}_{1},\ldots,\mathcal{D}_{N}\right)
g(x,t)  \right)  |_{x=0}\right\rangle _{\tau},\qquad g\in\mathcal{P}
_{V}.
\]
The form satisf\/ies $\left\langle f,g\right\rangle _{\tau}= \langle
wf,wg \rangle _{\tau}= \langle g,f \rangle _{\tau}$, $w\in W$.}

This is a general result for standard modules of the rational Cherednik
algebra, see~\cite{Dunkl/Opdam:2003}. The proof is based on induction on the
degree and the eigenfunction decomposition of the operator $\sum\limits_{i=1}
^{N}x_{i}\mathcal{D}_{i}$. Indeed
\begin{gather}
\sum_{i=1}^{N}x_{i}\mathcal{D}_{i}f(x,t)  = \langle
x,\nabla \rangle f(x,t)  +\sum_{v\in R_{+}}\kappa (
v )   \big(  f (  x,t\sigma_{v} )  -f (  x\sigma_{v}
,t\sigma_{v} )  \big)  , \label{sumxD}
\end{gather}
where $\nabla$ denotes the gradient (so that $ \langle x,\nabla
 \rangle =\sum\limits_{i=1}^{N}x_{i}\frac{\partial}{\partial x_{i}}$). Because
$\tau$ is irreducible there are integers $c_{\tau}(v)  $,
constant on conjugacy classes of ref\/lections (namely, the values of the
character of $\tau$) such that
\begin{gather*}
\sum_{v\in R_{+}}\kappa(v)  f (  x,t\sigma_{v} )
=\gamma (  \kappa;\tau )  f(x,t)  ,\qquad
\gamma (  \kappa;\tau )      :=\sum_{v\in R_{+}}c_{\tau} (
v )  \kappa(v),
\end{gather*}
for each $f\in\mathcal{P}_{V}$. The Laplacian is
\begin{gather*}
\Delta_{\kappa}f(x,t)      :=\sum_{i=1}^{N}\mathcal{D}_{i}%
^{2}f(x,t) \\
\hphantom{\Delta_{\kappa}f(x,t) }{}   =\Delta f(x,t)  +\sum_{v\in R_{+}}\kappa(v)
\left\{  2\frac{\langle v,\nabla f(  x,t\sigma_{v})
\rangle }{\langle x,v\rangle }-\vert v\vert
^{2}\frac{f(  x,t\sigma_{v})  -f(  x\sigma_{v},t\sigma
_{v} )  }{\langle x,v\rangle ^{2}}\right\}  ,
\end{gather*}
where $\Delta$ and $\nabla$ denote the ordinary Laplacian and gradient,
respectively. Motivated by the Gaussian inner product for scalar polynomials
(case $\tau=1$) which is def\/ined by
\[
 \langle f,g \rangle _{G}:=c_{\kappa}\int_{\mathbb{R}^{N}}f (
x )  g(x)  \prod_{v\in R_{+}} \vert  \langle
x,v \rangle  \vert ^{2\kappa(v)  }e^{- \vert
x \vert ^{2}/2}dx,
\]
where $c_{\kappa}$ is a normalizing (Macdonald--Mehta) constant, and satisf\/ies%
\[
 \langle f,g \rangle _{\tau}=\big\langle e^{-\Delta_{\kappa}
/2}f,e^{-\Delta_{\kappa}/2}g\big\rangle _{G},
\]
we def\/ine a bilinear \textit{Gaussian} form on $\mathcal{P}_{V}$ by
\[
 \langle f,g \rangle _{G}:=\big\langle e^{\Delta_{\kappa}
/2}f,e^{\Delta_{\kappa}/2}g\big\rangle _{\tau}.
\]
Note $e^{-\Delta_{\kappa}/2}:=\sum\limits_{n=0}^{\infty}\frac{(  -1/2)
^{n}}{n!}\Delta_{\kappa}^{n}$ is def\/ined for all polynomials since
$\Delta_{\kappa}$ is nilpotent. From the relations
\begin{gather*}
\Delta_{\kappa} (  x_{i}f(x,t)   )      =x_{i}
\Delta_{\kappa}f(x,t)  +2\mathcal{D}_{i}f(x,t)  ,\qquad
e^{-\Delta_{\kappa}/2} (  x_{i}f(x,t)  )      = (
x_{i}-\mathcal{D}_{i} )  e^{-\Delta_{\kappa}/2}f(x,t)  ,
\end{gather*}
we f\/ind that $ \langle  (  x_{i}-\mathcal{D}_{i} )
f,g \rangle _{G}= \langle f,\mathcal{D}_{i}g \rangle _{G}$, for
$1\leq i\leq N$ and $f,g\in\mathcal{P}_{V}$. Thus the multiplier operator
$x_{i}$ is self-adjoint for this form (since $x_{i}=\mathcal{D}_{i}
+\mathcal{D}_{i}^{\ast}$). This suggests that the form may have an expression
as an actual integral over~$\mathbb{R}^{N}$, at least for some restricted set
of the parameter values~$\kappa(v)  $. As in the scalar case
harmonic polynomials are involved in the analysis of the Gaussian form. The
equation
\[
\sum_{i=1}^{N} (  x_{i}\mathcal{D}_{i}+\mathcal{D}_{i}x_{i} )
=N+2 \langle x,\nabla \rangle +2\gamma (  \kappa;\tau )
\]
shows that
\[
\Delta_{\kappa} \vert x \vert ^{2m}f=2m \vert x \vert
^{2 (  m-1 )  } (  2m-2+N+2\gamma (  \kappa;\tau )
+2 \langle x,\nabla \rangle  )  f+ \vert x \vert
^{2m}\Delta_{\kappa}f,
\]
for $f\in\mathcal{P}_{V}$. For $n=0,1,2,\ldots$ let $\mathcal{P}%
_{V,n}= \{  f\in\mathcal{P}_{V}:f (  rx,t )  =r^{n}f (
x,t ),\forall\, r\in\mathbb{R} \}  $, the polynomials homogeneous of
degree $n$, and let $\mathcal{H}_{V,\kappa,n}:= \{  f\in\mathcal{P}
_{V,n}:\Delta_{\kappa}f=0 \}  $, the \textit{harmonic homogeneous}
polynomials. As a consequence of the previous formula, for $m=1,2,3,\ldots$
and $f\in\mathcal{H}_{V,\kappa,n}$ one obtains
\begin{gather}
\Delta_{\kappa}^{k}\big(   \vert x \vert ^{2m}f(x,t)
\big)  =4^{k} (  -m )  _{k} (  1-m-N/2-\gamma (
\kappa;\tau )  -n )  _{k} \vert x \vert ^{2m-2k}f (
x,t )  . \label{deltaF}
\end{gather}
(The Pochhammer symbol $(a)  _{k}$ is def\/ined by $(a)  _{0}\!=1$, $(a)  _{k+1}\!=(a)  _{k}(
a+k)  $ or \mbox{$(a)  _{k}:=\!\prod\limits_{i=1}^{k}(  a\!+\!i\!-\!1)
$}. In particular $\binom{n}{k}=(  -1)  ^{k}\frac{(  -n)
_{k}}{k!}$ and $(  -n)  _{k}=0$ for $n=0,\ldots,k-1$.) Thus
$\Delta_{\kappa}^{k}( \vert x\vert ^{2m}f(x,t)
 )  =0$ for $k>m$. With the same proofs as for the scalar case
\cite[Theorem~5.1.15]{Dunkl/Xu:2001} one can show (provided $\gamma (
\kappa;\tau )  +\frac{N}{2}\neq0,-1,-2,\ldots$): if $f\in\mathcal{P}
_{V,n}$ then
\begin{gather}
\pi_{\kappa,n}f     :=\sum_{j=0}^{ \lfloor n/2 \rfloor }\frac
{1}{4^{j}j! (  -\gamma (  \kappa;\tau )  -n+2-N/2 )  _{j}
} \vert x \vert ^{2j}\Delta_{\kappa}^{j}f\in\mathcal{H}_{V,\kappa
,n},\label{Hdecomp}\\
f     =\sum_{j=0}^{ \lfloor n/2 \rfloor }\frac{1}{4^{j}j! (
\gamma (  \kappa;\tau )  +N/2+n-2j )  _{j}} \vert
x \vert ^{2j}\pi_{\kappa,n-2j} \big(  \Delta_{\kappa}^{j}f \big)
.\nonumber
\end{gather}

From the def\/inition of $ \langle \cdot,\cdot \rangle _{\tau}$ it
follows that $f\in\mathcal{P}_{V,m}$, $g\in\mathcal{P}_{V,n}$ and $m\neq n$
implies $ \langle f,g \rangle _{\tau}=0$. Also if $f\in
\mathcal{H}_{V,\kappa,m}$, $g\in\mathcal{H}_{V,\kappa, n}$, $m\neq n$ then
$\big\langle  \vert x \vert ^{2a}f, \vert x \vert
^{2b}g\big\rangle _{\tau}=0$ for any $a,b=0,1,2,\ldots$. This follows from
the previous statement when $2a+m\neq2b+n$, otherwise $n=m+2a-2b$ and assume
$m<n$ (by symmetry of the form), thus $\big\langle  \vert x \vert
^{2a}f, \vert x \vert ^{2b}g\big\rangle _{\tau}=\big\langle
f,\Delta_{\kappa}^{a} \vert x \vert ^{2b}g\big\rangle _{\tau}=0$
because $a>b$. This shows that for generic $\kappa$ there is an orthogonal
decomposition of $\mathcal{P}_{V}$ as a~sum of $ \vert x \vert
^{2m}\mathcal{H}_{V,\kappa,n}$ over $m,n=0,1,2,\ldots$. If $f,g\in
\mathcal{H}_{V,\kappa,n}$ then
\begin{gather}\label{xsqform}
\big\langle  \vert x \vert ^{2m}f, \vert x \vert
^{2m}g\big\rangle _{\tau}     =\big\langle f,\Delta_{\kappa}^{m} \vert
x \vert ^{2m}g\big\rangle _{\tau}
   =4^{m}m!\left(  \frac{N}{2}+\gamma (  \kappa;\tau )  +n\right)
_{m} \langle f,g \rangle _{\tau},
\end{gather}
so to f\/ind an orthogonal basis for $\mathcal{P}_{V,m}$ it suf\/f\/ices to f\/ind an
orthogonal basis for each $\mathcal{H}_{V,\kappa,n}$.

The decomposition formula (\ref{Hdecomp}) implies the dimensionality result
for $\mathcal{H}_{V,\kappa,n}$ (when $\gamma (  \kappa;\tau )
+\frac{N}{2}\notin-\mathbb{N}_{0}$). It is clear that $\dim\mathcal{P}
_{V,n}=\dim\mathcal{P}_{\mathbb{R},n}\dim V=\binom{N+n-1}{N-1}\dim V$, and by
induction (for $n\geq1$)%
\begin{gather}
\dim\mathcal{H}_{V,\kappa,n}=\dim\mathcal{P}_{V,n}-\dim\mathcal{P}
_{V,n-2}=\binom{N+n-2}{N-2}\frac{N+2n-2}{N+n-2}\dim V. \label{dimH}
\end{gather}

We will need a lemma about integrating closed $1$-forms. Consider an $N$-tuple
$\mathbf{f}= (  f_{1},\ldots,f_{N} ) $ $\in \mathcal{P}_{V}^{N}$ as a
vector on which $W$ can act on the right. Say $\mathbf{f}$ is a \textit{closed
$1$-form} if $\mathcal{D}_{i}f_{j}-\mathcal{D}_{j}f_{i}=0$ for all~$i$,~$j$.

\begin{lemma}
Suppose $\mathbf{f}$ is a closed $1$-form and $1\leq j\leq N$ then
\[
\mathcal{D}_{j}\sum_{i=1}^{N}x_{i}f_{i}(x,t)  = (
 \langle x,\nabla \rangle +1+\gamma (  \kappa;\tau )
 )  f_{j}(x,t)  -\sum_{v\in R_{+}}\kappa(v)
 (  \mathbf{f} (  x\sigma_{v},t\sigma_{v} )  \sigma_{v} )
_{j}.
\]
\end{lemma}

\begin{proof}
By the commutation relations (\ref{dxxd})
\begin{gather*}
\mathcal{D}_{j}\sum_{i=1}^{N}x_{i}f_{i}(x,t)      =\sum_{i=1}
^{N}\left(  x_{i}\mathcal{D}_{j}f_{i}+\delta_{ij}f_{i}+2\sum_{v\in R_{+}
}\kappa(v)  \frac{v_{i}v_{j}}{ \vert v \vert ^{2}}
f_{i} (  x\sigma_{v},t\sigma_{v} )  \right) \\
\hphantom{\mathcal{D}_{j}\sum_{i=1}^{N}x_{i}f_{i}(x,t) }{}
  =\sum_{i=1}^{N}x_{i}\mathcal{D}_{i}f_{j}+f_{j}+2\sum_{v\in R_{+}}
\kappa(v)  \frac{v_{j}}{ \vert v \vert ^{2}}\sum
_{i=1}^{N}f_{i} (  x\sigma_{v},t\sigma_{v} )  v_{i}.
\end{gather*}
The calculation is f\/inished with the use of (\ref{sumxD}).
\end{proof}

\begin{corollary}\label{1-form}
Suppose $\mathbf{f}$ is a closed $1$-form,
homogeneous of degree $n$ and
\[
\sum\limits_{v\in R_{+}}\kappa (
v )  \mathbf{f} (  x\sigma_{v},t\sigma_{v} )  \sigma_{v}
=\lambda_{\kappa}\mathbf{f}(x,t)
\]
  for some constant
$\lambda_{\kappa}$ then
\[
\mathcal{D}_{j}\sum_{i=1}^{N}x_{i}f_{i}(x,t)  = (
n+1+\gamma (  \kappa;\tau )  -\lambda_{\kappa} )  f_{j} (
x,t )
\]
for $1\leq j\leq N$.
\end{corollary}

\section[The group $B_{2}$]{The group $\boldsymbol{B_{2}}$}\label{section3}

The rest of the paper concerns this group. The group $W (  B_{2} )  $
is the ref\/lection group with positive root system $R_{+}= \{   (
1,-1 )  , (  1,1 )  , (  1,0 )  , (  0,1 )
 \}  $. The corresponding ref\/lections are
\begin{gather*}
\sigma_{12}^{+}     :=\left[
\begin{matrix}
0 & 1\\
1 & 0
\end{matrix}
\right]  , \qquad \sigma_{12}^{-}:=\left[
\begin{matrix}
0 & -1\\
-1 & 0
\end{matrix}
\right]  ,\qquad
\sigma_{1}     :=\left[
\begin{matrix}
-1 & 0\\
0 & 1
\end{matrix}
\right]  , \qquad \sigma_{2}:=\left[
\begin{matrix}
1 & 0\\
0 & -1
\end{matrix}
\right]  .
\end{gather*}
The values of $\kappa$ on the conjugacy classes $ \{  \sigma_{12}%
^{+},\sigma_{12}^{-} \}  $ and $ \{  \sigma_{1},\sigma_{2} \}  $
will be denoted by~$k_{0}$ and~$k_{1}$, respectively. We consider the unique
$2$-dimensional representation $\tau$ and set $V:=\mathrm{span} \{
t_{1},t_{2} \}  $. The generic element of $\mathcal{P}_{V}$ is $f (
x,t )  =f_{1}(x)  t_{1}+f_{2}(x)  t_{2}$. The
ref\/lections act on this polynomial as follows
\begin{gather*}
\sigma_{12}^{+}f(x,t)      =f_{2} (  x_{2},x_{1} )
t_{1}+f_{1} (  x_{2},x_{1} )  t_{2},\\
\sigma_{12}^{-}f(x,t)     =-f_{2} (  -x_{2},-x_{1} )
t_{1}-f_{1} (  -x_{2},-x_{1} )  t_{2},\\
\sigma_{1}f(x,t)      =-f_{1} (  -x_{1},x_{2} )
t_{1}+f_{2} (  -x_{1},x_{2} )  t_{2},\\
\sigma_{2}f(x,t)      =f_{1} (  x_{1},-x_{2} )
t_{1}-f_{2} (  x_{1},-x_{2} )  t_{2}.
\end{gather*}
Here is the formula for $\mathcal{D}_{1}$ ($\mathcal{D}_{2}$ is similar)
\begin{gather*}
\mathcal{D}_{1}f(x,t)  =\frac{\partial}{\partial x_{1}}f (
x,t )  +k_{1}\frac{f (  x_{1},x_{2},-t_{1},t_{2} )  -f (
-x_{1},x_{2},-t_{1},t_{2} )  }{x_{1}}\\
\hphantom{\mathcal{D}_{1}f(x,t)  =}{}
+k_{0}\frac{f (  x_{1},x_{2},t_{2},t_{1} )  -f (  x_{2}
,x_{1},t_{2},t_{1} )  }{x_{1}-x_{2}}\\
\hphantom{\mathcal{D}_{1}f(x,t)  =}{}
+k_{0}\frac{f (  x_{1},x_{2},-t_{2},-t_{1} )  -f (  -x_{2},-x_{1},-t_{2},-t_{1} )  }{x_{1}+x_{2}}.
\end{gather*}
Since the matrices for the ref\/lections all have trace zero we f\/ind that
$\gamma (  \kappa;\tau )  =0$. We investigate the properties of
$ \langle \cdot,\cdot \rangle _{\tau}$ by constructing bases for each
$\mathcal{H}_{V,\kappa,n}$; note $\dim\mathcal{H}_{V,\kappa,n}=4$ for $n\geq
1$. The form on $V$ is given by $ \langle t_{i},t_{j} \rangle _{\tau
}=\delta_{ij}$. A convenient orthogonal decomposition is based on
\textit{types} def\/ined by the action of $ (  \sigma_{1},\sigma_{2} )
$: suppose $f$ satisf\/ies $\sigma_{1}f=\varepsilon_{1}f$, $\sigma_{2}
f=\varepsilon_{2}f$, then $f$ is of type EE, EO, OE, OO if $ (
\varepsilon_{1},\varepsilon_{2} )  = (  1,1 )  , (
1,-1 )  , (  -1,1 )  , (  -1,-1 )  $ respectively (E~is
for even, O~is for odd). Because the ref\/lections are self-adjoint for the form
$ \langle \cdot,\cdot \rangle _{\tau}$ it follows immediately that
polynomials of dif\/ferent types are mutually orthogonal. By formula~(\ref{dimH}) $\dim\mathcal{H}_{V,\kappa,n}=4$ for $n\geq1$. The construction
of the harmonic basis $ \{  p_{n,i}:n\geq0,\,1\leq i\leq4 \}  $ is by a
recurrence; the process starts at degree~$1$. Set
\begin{gather*}
p_{1,1}     :=x_{1}t_{1}+x_{2}t_{2},\qquad
p_{1,2}     :=-x_{2}t_{1}+x_{1}t_{2},\\
p_{1,3}     :=x_{1}t_{1}-x_{2}t_{2},\qquad
p_{1,4}     :=-x_{2}t_{1}-x_{1}t_{2}.
\end{gather*}

Thus $p_{1,1}$, $p_{1,3}$ are of type EE and $p_{1,2}$, $p_{1,4}$ are of type OO.
Also $ \langle p_{1,1},p_{1,3} \rangle _{\tau}=0= \langle
p_{1,2},p_{1,4} \rangle _{\tau}$ because $\sigma_{12}^{+}p_{1,1}=p_{1,1}
$, $\sigma_{12}^{+}p_{1,3}=-p_{1,3}$, $\sigma_{12}^{+}p_{1,2}=-p_{1,2}$,
$\sigma_{12}^{+}p_{1,4}=p_{1,4}$. The same decompositions work in
$\mathcal{H}_{V,\kappa,n}$ for each odd~$n$. By direct computation we obtain%
\begin{gather*}
 \langle p_{1,1},p_{1,1} \rangle _{\tau}     =2 (  1-2k_{0}-2k_{1} )  ,\qquad
 \langle p_{1,2},p_{1,2} \rangle _{\tau}     =2 (  1+2k_{0}+2k_{1} )  ,\\
 \langle p_{1,3},p_{1,3} \rangle _{\tau}     =2 (  1+2k_{0}-2k_{1} )  ,\qquad
 \langle p_{1,4},p_{1,4} \rangle _{\tau}     =2 (  1-2k_{0}+2k_{1} )  .
\end{gather*}
These values show that a necessary condition for the existence of a
positive-def\/inite inner product on $\mathcal{P}_{V}$ is $-\frac{1}{2}<\pm
k_{0}\pm k_{1}<\frac{1}{2}$. Note this is signif\/icantly dif\/ferent from the
analogous condition for scalar polynomials (that is, $\tau=1$), namely
$k_{0},k_{1},k_{0}+k_{1}>-\frac{1}{2}$. The formulae become more readable by
the use of
\[
k_{+}:=k_{0}+k_{1},\qquad k_{-}:=k_{1}-k_{0}.
\]
The recurrence starts at degree $0$ by setting $p_{0,1}=p_{0,3}=t_{1}$ and
$p_{0,2}=-p_{0,4}=t_{2}$ (in the exceptional case $n=0$ there are only two
linearly independent polynomials).

\begin{definition}
\label{pdef}The polynomials $p_{n,i}\in\mathcal{P}_{V,n}$ for $n\geq1$ and
$1\leq i\leq4$ are given by%
\[
\left[
\begin{matrix}
p_{2m+1,1} & p_{2m+1,3}\\
p_{2m+1,2} & p_{2m+1,4}
\end{matrix}
\right]  :=\left[
\begin{matrix}
x_{1} & x_{2}\\
-x_{2} & x_{1}%
\end{matrix}
\right]  \left[
\begin{matrix}
p_{2m,1} & p_{2m,3}\\
p_{2m,2} & p_{2m,4}
\end{matrix}
\right]  , \qquad m\geq0,
\]
and%
\begin{gather*}
\left[
\begin{matrix}
p_{2m,1}\\
p_{2m,2}
\end{matrix}
\right]     :=\frac{1}{2m-1}\left[
\begin{matrix}
x_{1} & x_{2}\\
-x_{2} & x_{1}
\end{matrix}
\right]  \left[
\begin{matrix}
(  2m-1+2k_{-})  p_{2m-1,3}\\
(  2m-1-2k_{-})  p_{2m-1,4}
\end{matrix}
\right]  ,\\
\left[
\begin{matrix}
p_{2m,3}\\
p_{2m,4}
\end{matrix}
\right]     :=\frac{1}{2m-1}\left[
\begin{matrix}
x_{1} & x_{2}\\
-x_{2} & x_{1}
\end{matrix}
\right]  \left[
\begin{matrix}
(  2m-1+2k_{+})  p_{2m-1,1}\\
(  2m-1-2k_{+})  p_{2m-1,2}
\end{matrix}
\right]  ,
\end{gather*}
for $m\geq1$.
\end{definition}

The construction consists of two disjoint sequences
\begin{gather*}
\left[
\begin{matrix}
p_{1,1}\\
p_{1,2}%
\end{matrix}
\right]  \overset{1\pm2k_{+}}{\rightarrow}\left[
\begin{matrix}
p_{2,3}\\
p_{2,4}%
\end{matrix}
\right]      \rightarrow\left[
\begin{matrix}
p_{3,3}\\
p_{3,4}
\end{matrix}
\right]  \overset{3\pm2k_{-}}{\rightarrow}\left[
\begin{matrix}
p_{4,1}\\
p_{4,2}
\end{matrix}
\right]  \rightarrow\left[
\begin{matrix}
p_{5,1}\\
p_{5,2}
\end{matrix}
\right]  \overset{5\pm2k_{+}}{\rightarrow},\\
\left[
\begin{matrix}
p_{1,3}\\
p_{1,4}
\end{matrix}
\right]  \overset{1\pm2k_{-}}{\rightarrow}\left[
\begin{matrix}
p_{2,1}\\
p_{2,2}
\end{matrix}
\right]     \rightarrow\left[
\begin{matrix}
p_{3,1}\\
p_{3,2}
\end{matrix}
\right]  \overset{3\pm2k_{+}}{\rightarrow}\left[
\begin{matrix}
p_{4,3}\\
p_{4,4}
\end{matrix}
\right]  \rightarrow\left[
\begin{matrix}
p_{5,3}\\
p_{5,4}
\end{matrix}
\right]  \overset{5\pm2k_{-}}{\rightarrow}.
\end{gather*}

The typical steps are $ [  p_{2m+1,\ast} ]  =X [  p_{2m,\ast
} ]  $ and $ [  p_{2m+2,\ast} ]  =XC_{2m+1} [
p_{2m+1,\circ} ]$ ($\ast$ and $\circ$ denote a $2$-vector with labels
$1$, $2$ or $3$, $4$ as in the diagram) where
\[
X=\left[
\begin{matrix}
x_{1} & x_{2}\\
-x_{2} & x_{1}
\end{matrix}
\right]  , \qquad C_{2m+1}=\left[
\begin{matrix}
\frac{2m+1+2\lambda}{2m+1} & 0\\
0 & \frac{2m+1-2\lambda}{2m+1}
\end{matrix}
\right]  ,
\]
where $\lambda=k_{+}$ or $k_{-}$ (indicated by the labels on the arrows). Each
polynomial is nonzero, regardless of the parameter values. This follows from
evaluation at $x=(  1,\mathrm{i})$. Let $z:=t_{1}+\mathrm{i}t_{2}$.

\begin{proposition}\label{pnonz}
If $n=4m+1$ or $4m$ with $m\geq1$ then
\[
\left[  p_{n,1}\left(  1,\mathrm{i}\right)  ,p_{n,2}\left(  1,\mathrm{i}%
\right)  ,p_{n,3}\left(  1,\mathrm{i}\right)  ,p_{n,4}\left(  1,\mathrm{i}%
\right)  \right]  =2^{n-1}\left[  z,-\mathrm{i}z,\overline{z},-\mathrm{i}%
\overline{z}\right]  ;
\]
if $n=4m+2$ or $4m+3$ with $m\geq0$ then
\[
\left[  p_{n,1}\left(  1,\mathrm{i}\right)  ,p_{n,2}\left(  1,\mathrm{i}%
\right)  ,p_{n,3}\left(  1,\mathrm{i}\right)  ,p_{n,4}\left(  1,\mathrm{i}%
\right)  \right]  =2^{n-1}\left[  \overline{z},-\mathrm{i}\overline
{z},z,-\mathrm{i}z\right]  .
\]

\end{proposition}

\begin{proof}
The formula is clearly valid for $n=1$. The typical step in the inductive
proof is
\[
\begin{bmatrix}
p_{n+1,1}\\
p_{n+1,2}%
\end{bmatrix}
=\left[
\begin{matrix}
1 & \mathrm{i}\\
-\mathrm{i} & 1
\end{matrix}
\right]
\begin{bmatrix}
c_{1}p_{n,\ell}\\
-\mathrm{i}c_{2}p_{n,\ell}
\end{bmatrix}
=2
\begin{bmatrix}
p_{n,\ell}\\
-\mathrm{i}p_{n,\ell}
\end{bmatrix}
,
\]
because $c_{1}+c_{2}=2$. If $n$ is odd then $\ell=3$, otherwise $\ell=1$. The
same argument works for $p_{n+1,3}$, $p_{n+1,4}$.
\end{proof}

The types of $p_{2m+1,1}$, $p_{2m+1,2}$, $p_{2m+1,3}$, $p_{2m+1,4}$ are EE, OO, EE, OO
respectively and the types of $p_{2m,1}$, $p_{2m,2}$, $p_{2m,3}$, $p_{2m,4}$ are OE,
EO, OE, EO respectively.

\begin{proposition}\label{actps12}
For $m\geq0$
\[
\sigma_{12}^{+}\left[
\begin{matrix}%
p_{2m+1,1}\\
p_{2m+1,2}\\
p_{2m+1,3}\\
p_{2m+1,4}
\end{matrix}
\right]  =\left[
\begin{matrix}
p_{2m+1,1}\\
-p_{2m+1,2}\\
-p_{2m+1,3}\\
p_{2m+1,4}
\end{matrix}
\right]  ,\qquad \sigma_{12}^{+}\left[
\begin{matrix}%
p_{2m,1}\\
p_{2m,2}\\
p_{2m,3}\\
p_{2m,4}%
\end{matrix}
\right]  =\left[
\begin{matrix}
p_{2m,2}\\
p_{2m,1}\\
-p_{2m,4}\\
-p_{2m,3}
\end{matrix}
\right]  .
\]
\end{proposition}

\begin{proof}
Using induction one needs to show that the validity of the statements for
$2m-1$ implies the validity for $2m$, and the validity for $2m$ implies the
validity for $2m+1$, for each $m\geq1$. Suppose the statements hold for some~$2m$. The types of $p_{2m+1,i}$ are easy to verify ($1\leq i\leq4$). Next
\begin{gather*}
\sigma_{12}^{+}p_{2m+1,1}     =\sigma_{12}^{+} (  x_{1}p_{2m,1}%
+x_{2}p_{2m,2} )  =x_{2}\sigma_{12}^{+}p_{2m,1}+x_{1}\sigma_{12}%
^{+}p_{2m,2}\\
\hphantom{\sigma_{12}^{+}p_{2m+1,1}}{}
 =x_{2}p_{2m,2}+x_{1}p_{2m,1}=p_{2m+1,1},\\
\sigma_{12}^{+}p_{2m+1,2}     =\sigma_{12}^{+} (  x_{1}p_{2m,2}%
-x_{2}p_{2m,1} )  =x_{2}\sigma_{12}^{+}p_{2m,2}-x_{1}\sigma_{12}%
^{+}p_{2m,1}
   =-p_{2m+1,2},
\end{gather*}
and by similar calculations $\sigma_{12}^{+}p_{2m+1,3}=-p_{2m+1,3}$ and
$\sigma_{12}^{+}p_{2m+1,4}=-p_{2m+1,4}$. Now suppose the statements hold for
some $2m-1$. As before the types of $p_{2m,i}$ are easy to verify. Consider
\begin{gather*}
\sigma_{12}^{+}p_{2m,1}     =\sigma_{12}^{+}\left(  \frac{2m-1+2k_{-}}
{2m-1}x_{1}p_{2m-1,3}+\frac{2m-1-2k_{-}}{2m-1}x_{2}p_{2m-1,4}\right) \\
\hphantom{\sigma_{12}^{+}p_{2m,1}}{}
  =-\frac{2m-1+2k_{-}}{2m-1}x_{2}p_{2m-1,3}+\frac{2m-1-2k_{-}}{2m-1}%
x_{1}p_{2m-1,4}   =p_{2m,2},\\
\sigma_{12}^{+}p_{2m,3}     =\sigma_{12}^{+}\left(  \frac{2m-1+2k_{+}}
{2m-1}x_{1}p_{2m-1,1}+\frac{2m-1-2k_{+}}{2m-1}x_{2}p_{2m-1,2}\right) \\
\hphantom{\sigma_{12}^{+}p_{2m,3}}{}
 =\frac{2m-1+2k_{+}}{2m-1}x_{2}p_{2m-1,1}-\frac{2m-1-2k_{+}}{2m-1}%
x_{1}p_{2m-1,2}
   =-p_{2m,4}.
\end{gather*}
Since $ (  \sigma_{12}^{+} )  ^{2}=1$ this f\/inishes the inductive proof.
\end{proof}

\begin{corollary}
The polynomials $ \{  p_{2m+1,i}:1\leq i\leq4 \}  $ are mutually orthogonal.
\end{corollary}

\begin{proof}
These polynomials are eigenfunctions of the (self-adjoint) ref\/lections
$\sigma_{1}$, $\sigma_{12}^{+}$ with dif\/fe\-rent pairs of eigenvalues.
\end{proof}

We will use Corollary \ref{1-form} to evaluate $\mathcal{D}_{i}p_{n,j}$ (with
$i=1,2$ and $1\leq j\leq4$). We intend to show that the expressions appearing
in Def\/inition~\ref{pdef} can be interpreted as closed $1$-forms $\mathbf{f}$
which are eigenfunctions of $\mathbf{f\longmapsto}\sum\limits_{v\in R_{+}}
\kappa(v)  \mathbf{f} (  x\sigma_{v},t\sigma_{v} )
\sigma_{v}$. In the present context this is the operator
\begin{gather*}
\rho (  f_{1}(x,t)  ,f_{2}(x,t)  )
:=k_{1} (   (  -\sigma_{1}+\sigma_{2} )  f_{1}, (  \sigma
_{1}-\sigma_{2} )  f_{2} )
   +k_{0} (   (  \sigma_{12}^{+}-\sigma_{12}^{-} )  f_{2}, (
\sigma_{12}^{+}-\sigma_{12}^{-} )  f_{1} )  .
\end{gather*}
It is easy to check that $ (  \sigma_{1}-\sigma_{2} )  f (
x,t )  =0= (  \sigma_{12}^{+}-\sigma_{12}^{-} )  f (
x,t )  $ whenever $f\in\mathcal{P}_{V,2m+1}$ for $m\geq0$. For the even
degree polynomials $f(x,t)  $ one f\/inds $ (  \sigma_{12}%
^{+}-\sigma_{12}^{-} )  f=2\sigma_{12}^{+}f$ , $ (  \sigma_{1}
-\sigma_{2} )  f=-2f$ if $f$ is of type OE and $ (  \sigma_{1}-\sigma_{2} )  f=2f$ if $f$ is of type EO.

\begin{theorem}
For $m\geq1$
\begin{gather*}
\mathcal{D}_{1}p_{2m,1}     =2m\frac{2m-1+2k_{-}}{2m-1}p_{2m-1,3}
,\qquad \mathcal{D}_{2}p_{2m,1}=2m\frac{2m-1-2k_{-}}{2m-1}p_{2m-1,4},\\
\mathcal{D}_{1}p_{2m,2}     =2m\frac{2m-1-2k_{-}}{2m-1}p_{2m-1,4}
,\qquad \mathcal{D}_{2}p_{2m,2}=-2m\frac{2m-1+2k_{-}}{2m-1}p_{2m-1,3},\\
\mathcal{D}_{1}p_{2m,3}     =2m\frac{2m-1+2k_{+}}{2m-1}p_{2m-1,1}
,\qquad \mathcal{D}_{2}p_{2m,3}=2m\frac{2m-1-2k_{+}}{2m-1}p_{2m-1,2},\\
\mathcal{D}_{1}p_{2m,4}     =2m\frac{2m-1-2k_{+}}{2m-1}p_{2m-1,2}
,\qquad \mathcal{D}_{2}p_{2m,4}=-2m\frac{2m-1+2k_{+}}{2m-1}p_{2m-1,1},
\end{gather*}
and for $m\geq1$ $($in vector notation$)$
\begin{gather*}
 [  \mathcal{D}_{1},\mathcal{D}_{2} ]  p_{2m+1,1}     = (
2m+1-2k_{+} )   [  p_{2m,1},p_{2m,2} ]  ,\\
 [  \mathcal{D}_{1},\mathcal{D}_{2} ]  p_{2m+1,2}     = (
2m+1+2k_{+} )   [  p_{2m,2},-p_{2m,1} ]  ,\\
 [  \mathcal{D}_{1},\mathcal{D}_{2} ]  p_{2m+1,3}     = (
2m+1-2k_{-} )   [  p_{2m,3},p_{2m,4} ]  ,\\
 [  \mathcal{D}_{1},\mathcal{D}_{2} ]  p_{2m+1,4}     = (
2m+1+2k_{-} )   [  p_{2m,4},-p_{2m,3} ]  .
\end{gather*}
\end{theorem}

\begin{proof}
Use induction as above. In each case write $p_{n,i}$ in the form $x_{1}
f_{1}+x_{2}f_{2}$ with $\mathcal{D}_{2}f_{1}=\mathcal{D}_{1}f_{1}$ and apply
Corollary~\ref{1-form}. Suppose the statements are true for some $2m-1$
(direct verif\/ication for $m=1$). Then
\[
p_{2m,1}=x_{1}\left(  \frac{2m-1+2k_{-}}{2m-1}p_{2m-1,3}\right)  +x_{2}\left(
\frac{2m-1-2k_{-}}{2m-1}p_{2m-1,4}\right)
\]
and
\begin{gather*}
\mathcal{D}_{2}\left(  \frac{2m-1+2k_{-}}{2m-1}p_{2m-1,3}\right)
=\frac{\left(  2m-1+2k_{-}\right)  \left(  2m-1-2k_{-}\right)  }
{2m-1}p_{2m-2,4}\\
 \hphantom{\mathcal{D}_{2}\left(  \frac{2m-1+2k_{-}}{2m-1}p_{2m-1,3}\right)  }{}
   =\mathcal{D}_{1}\left(  \frac{2m-1-2k_{-}}{2m-1}p_{2m-1,4}\right)
\end{gather*}
by the inductive hypothesis. Since $\rho p_{2m-1,i}=0$ for $1\leq i\leq4$
Corollary \ref{1-form} shows that $\mathcal{D}_{1}p_{2m,1}=2m\big(
\frac{2m-1+2k_{-}}{2m-1}p_{2m-1,3}\big)  $ and $\mathcal{D}_{2}p_{2m,1}$ has
the stated value. Also
\[
p_{2m,1}=x_{1}\left(  \frac{2m-1-2k_{+}}{2m-1}p_{2m-1,2}\right)  +x_{2}\left(
-\frac{2m-1+2k_{+}}{2m-1}p_{2m-1,1}\right)  ,
\]
and
\begin{gather*}
\mathcal{D}_{2}\left(  \frac{2m-1-2k_{+}}{2m-1}p_{2m-1,2}\right)
=-\frac{\left(  2m-1+2k_{-}\right)  \left(  2m-1-2k_{+}\right)  }
{2m-1}p_{2m-2,1}\\
\hphantom{\mathcal{D}_{2}\left(  \frac{2m-1-2k_{+}}{2m-1}p_{2m-1,2}\right) }{}
  =\mathcal{D}_{1}\left(  -\frac{2m-1+2k_{+}}{2m-1}p_{2m-1,2}\right)  ,
\end{gather*}
and so $\mathcal{D}_{2}p_{2m,4}=2m\big(  {-}\frac{2m-1+2k_{+}}{2m-1}%
p_{2m-1,1}\big)  $. The other statements for $p_{2m,i}$ have similar proofs.
Now assume the statements are true for some $2m$. The desired relations follow
from Corollary~\ref{1-form} and the following
\begin{gather*}
\rho (  p_{2m,1},p_{2m,2} )      =2k_{0} (  \sigma_{12}%
^{+}p_{2m,2},\sigma_{12}^{+}p_{2m,1} )  +2k_{1} (  p_{2m,1}
,p_{2m,2} ) \\
\hphantom{\rho (  p_{2m,1},p_{2m,2} ) }{}
 = (  2k_{0}+2k_{1} )   (  p_{2m,1},p_{2m,2} )
,\qquad \lambda_{\kappa}=2k_{+};\\
\rho (  p_{2m,2},-p_{2m,1} )     =2k_{0} (  -\sigma_{12}
^{+}p_{2m,1},\sigma_{12}^{+}p_{2m,2} )  -2k_{1} (  p_{2m,2}
,-p_{2m,1} ) \\
\hphantom{\rho (  p_{2m,2},-p_{2m,1} )}{}
  = (  -2k_{0}-2k_{1} )  (  p_{2m,2},-p_{2m,1} )
,\qquad \lambda_{\kappa}=-2k_{+};\\
\rho (  p_{2m,3},p_{2m,4} )      =2k_{0} (  \sigma_{12}%
^{+}p_{2m,4},\sigma_{12}^{+}p_{2m,3} )  +2k_{1} (  p_{2m,3}
,p_{2m,4} ) \\
\hphantom{\rho (  p_{2m,3},p_{2m,4} )}{}
   = (  -2k_{0}+2k_{1} )   (  p_{2m,3},p_{2m,4} )
,\qquad \lambda_{\kappa}=2k_{-};\\
\rho (  p_{2m,4},-p_{2m,3} )      =2k_{0} (  -\sigma_{12}%
^{+}p_{2m,3},\sigma_{12}^{+}p_{2m,4} )  -2k_{1} (  p_{2m,4}
,-p_{2m,3} ) \\
\hphantom{\rho (  p_{2m,4},-p_{2m,3} )}{}
   = (  2k_{0}-2k_{1} )  (  p_{2m,4},-p_{2m,3} )
,\qquad \lambda_{\kappa}=-2k_{-}.
\end{gather*}
The closed $1$-form condition is shown by using the inductive hypothesis. A
typical calculation is
\begin{gather*}
p_{2m+1,2}     =x_{1}p_{2m,2}-x_{2}p_{2m,1},\qquad
\mathcal{D}_{2}p_{2m,2}     =-2m\frac{2m-1+2k_{-}}{2m-1}p_{2m-1,3}
=-\mathcal{D}_{1}p_{2m,1}.
\end{gather*}
This completes the inductive proof.
\end{proof}

\begin{corollary}
Suppose $n\geq1$ and $1\leq i\leq4$ then $\Delta_{\kappa}p_{n,i}=0$.
\end{corollary}

\begin{proof}
By the above formulae $\mathcal{D}_{1}^{2}p_{n,i}=-\mathcal{D}_{2}^{2}p_{n,i}$
in each case. As a typical calculation%
\begin{gather*}
\mathcal{D}_{1}^{2}p_{2m,1}     =\mathcal{D}_{1}\left(  2m\frac{2m-1+2k_{-}
}{2m-1}p_{2m-1,3}\right) \\
\hphantom{\mathcal{D}_{1}^{2}p_{2m,1} }{}
  =\frac{2m}{2m-1} (  2m-1+2k_{-} )   (  2m+1-2k_{-} )
p_{2m-2,3},\\
\mathcal{D}_{2}^{2}p_{2m,1}     =\mathcal{D}_{2}\left(  2m\frac{2m-1-2k_{-}%
}{2m-1}p_{2m-1,4}\right) \\
\hphantom{\mathcal{D}_{2}^{2}p_{2m,1} }{}
  =\frac{2m}{2m-1} (  2m+1-2k_{-} )   (  2m+1+2k_{-} )
 (  -p_{2m-2,3} )  .
\end{gather*}
It suf\/f\/ices to check the even degree cases because $\Delta_{\kappa}f=0$
implies $\Delta_{\kappa}\mathcal{D}_{1}f=0$.
\end{proof}

Because $p_{2m,1}$ and $p_{2m,3}$ are both type OE one can use an appropriate
self-adjoint operator to prove orthogonality. Indeed let $\mathcal{U}%
_{12}:=\sigma_{12}^{+} (  x_{2}\mathcal{D}_{1}-x_{1}\mathcal{D}%
_{2} )  $, which is self-adjoint for the form $ \langle \cdot
,\cdot \rangle _{\tau}$.

\begin{proposition}\label{actu12p}
Suppose $m\geq1$ and $1\leq i\leq4$ then $\mathcal{U}%
_{12}p_{2m,i}=2m\varepsilon_{i}p_{2m,i}$ where $\varepsilon_{1}=\varepsilon
_{4}=-1$ and $\varepsilon_{2}=\varepsilon_{3}=1$.
\end{proposition}

\begin{proof}
This is a simple verif\/ication. For example
\begin{gather*}
 (  x_{2}\mathcal{D}_{1}-x_{1}\mathcal{D}_{2} )  p_{2m,4}
=\frac{2m}{2m-1} (  x_{2} (  2m-1-2k_{+} )  p_{2m-1,3}
+x_{1} (  2m-1+2k_{+} )  p_{2m-1,4} ) \\
 \phantom{(  x_{2}\mathcal{D}_{1}-x_{1}\mathcal{D}_{2} )  p_{2m,4}}{}
  =2mp_{2m,3},
\end{gather*}
and $\sigma_{12}^{+}p_{2m,3}=-p_{2m,4}$.
\end{proof}

\begin{corollary}
For $m\geq1$ the polynomials $ \{  p_{2m,i}:1\leq i\leq4 \}  $ are
mutually orthogonal.
\end{corollary}

\begin{proof}
These polynomials are eigenfunctions of the self-adjoint operators $\sigma
_{1}$ and $\mathcal{U}_{12}$ with dif\/ferent pairs of eigenvalues.
\end{proof}

\begin{proposition}
\label{Hbasis}For any $k_{0}$, $k_{1}$ and $n\geq1$ the polynomials $ \{
p_{n,i}:1\leq i\leq4 \}  $ form a basis for~$\mathcal{H}_{V,\kappa,n}$.
\end{proposition}

\begin{proof}
The fact that $p_{n,i}\neq0$ for all $n$, $i$ (Proposition~\ref{pnonz}) and the
eigenvector properties from Propositions~\ref{actps12} and~\ref{actu12p} imply
linear independence for each set $\{  p_{n,i}:1\leq i\leq4\}$.
\end{proof}

We introduce the notation $\nu(f)  := \langle
f,f \rangle _{\tau}$ (without implying that the form is positive). To
evaluate $\nu (  p_{n,i} )  $ we will use induction and the following
simple fact: suppose $f=x_{1}f_{1}+x_{2}f_{2}$ and $\mathcal{D}_{i}f=\alpha
f_{i}$ for $i=1,2$ and some $\alpha\in\mathbb{Q} [  k_{0},k_{1} ]  $
then
\begin{gather*}
\nu(f)      = \langle x_{1}f_{1},f \rangle _{\tau
}+ \langle x_{2}f_{2},f \rangle _{\tau}= \langle f_{1}
,\mathcal{D}_{1}f \rangle _{\tau}+ \langle f_{2},\mathcal{D}
_{2}f \rangle _{\tau}   =\alpha (  \nu (  f_{1} )  +\nu (  f_{2} )   ).
\end{gather*}
By use of the various formulae we obtain
\begin{gather*}
\nu (  p_{2m,1} )      =\nu (  p_{2m,2} )
   =2m\!\left(  \left(  \frac{2m\!-\!1\!+\!2k_{-}}{2m\!-\!1}\right)  ^{2}\!\nu (
p_{2m-1,3} )  +\!\left(  \frac{2m\!-\!1\!-\!2k_{-}}{2m\!-\!1}\right)  ^{2}\!\nu (
p_{2m-1,4} )  \right) \! ,\\
\nu (  p_{2m,3} )      =\nu (  p_{2m,4} )
   =2m\!\left(  \left(  \frac{2m\!-\!1\!+\!2k_{+}}{2m\!-\!1}\right)  ^{2}\!\nu (
p_{2m-1,1} )  +\!\left(  \frac{2m\!-\!1\!-\!2k_{+}}{2m\!-\!1}\right)  ^{2}\!\nu (
p_{2m-1,2} )  \right)\!,
\end{gather*}
and%
\begin{gather*}
\nu (  p_{2m+1,1} )      =2 (  2m+1-2k_{+} )  \nu (p_{2m,1} )  ,\qquad
\nu (  p_{2m+1,2} )      =2 (  2m+1+2k_{+} )  \nu (
p_{2m,1} )  ,\\
\nu (  p_{2m+1,3} )      =2 (  2m+1-2k_{-} )  \nu (
p_{2m,3} )  ,\qquad
\nu (  p_{2m+1,4} )      =2 (  2m+1+2k_{-} )  \nu (
p_{2m,3} )  .
\end{gather*}
From these relations it follows that
\begin{gather*}
\nu (  p_{2m+2,1} )      =8 (  m+1 )  \frac{(
2m+1-2k_{-} )   (  2m+1+2k_{-} )  }{2m+1}\nu (
p_{2m,3} )  ,\\
\nu (  p_{2m+2,3} )      =8 (  m+1 )  \frac{ (
2m+1-2k_{+} )   (  2m+1+2k_{+} )  }{2m+1}\nu (
p_{2m,1} )  .
\end{gather*}
It is now a case-by-case verif\/ication for an inductive proof. Def\/ine a
Pochhammer-type function (referred to as ``$\Pi$'' in the sequel) for notational convenience
\[
\Pi\left(  a,b,m,\varepsilon_{1}\varepsilon_{2}\varepsilon_{3}\varepsilon
_{4}\right)  :=\frac{\left(  \frac{1}{4}+\frac{a}{2}\right)  _{m+\varepsilon
_{1}}\left(  \frac{1}{4}-\frac{a}{2}\right)  _{m+\varepsilon_{2}}\left(
\frac{3}{4}+\frac{b}{2}\right)  _{m+\varepsilon_{3}}\left(  \frac{3}{4}%
-\frac{b}{2}\right)  _{m+\varepsilon_{4}}}{\left(  \frac{1}{4}\right)
_{m+\varepsilon_{1}}\left(  \frac{1}{4}\right)  _{m+\varepsilon_{2}}\left(
\frac{3}{4}\right)  _{m+\varepsilon_{3}}\left(  \frac{3}{4}\right)
_{m+\varepsilon_{4}}},
\]
where $m\geq0$ and $\varepsilon_{1}\varepsilon_{2}\varepsilon_{3}%
\varepsilon_{4}$ is a list with $\varepsilon_{i}=0$ or $1$. To avoid
repetitions of the factor $2^{n}n!$ let
\[
\nu^{\prime} (  p_{n,i} )  :=\nu (  p_{n,i} )  / (
2^{n}n! )  .
\]
Then
\begin{gather}
\nu^{\prime} (  p_{4m,1} )      =\nu^{\prime} (  p_{4m,2} )
=\Pi (  k_{+},k_{-},m,0000 )  ,\nonumber\\
\nu^{\prime} (  p_{4m,3} )      =\nu^{\prime} (  p_{4m,4} )
=\Pi (  k_{-},k_{+},m,0000 )  ,\nonumber\\
\nu^{\prime} (  p_{4m+2,1} )      =\nu^{\prime} (  p_{4m+2,2}  )  =\Pi (  k_{-},k_{+},m,1100 )  ,\nonumber\\
\nu^{\prime} (  p_{4m+2,3} )      =\nu^{\prime} (  p_{4m+2,4}  )  =\Pi (  k_{+},k_{-},m,1100 )  ,\label{norms1}
\end{gather}
and
\begin{gather}
\nu^{\prime} (  p_{4m+1,1} )      =\Pi (  k_{+},k_{-}
,m,0100 )  ,\qquad \nu^{\prime} (  p_{4m+1,2} )  =\Pi (
k_{+},k_{-},m,1000 )  ,\nonumber\\
\nu^{\prime} (  p_{4m+1,3} )      =\Pi (  k_{-},k_{+}
,m,0100 )  ,\qquad \nu^{\prime} (  p_{4m+1,4} )  =\Pi (
k_{-},k_{+},m,1000 )  ,\label{nrm2}
\\
\nu^{\prime} (  p_{4m+3,1} )      =\Pi (  k_{-},k_{+}
,m,1101 )  ,\qquad \nu^{\prime} (  p_{4m+3,2} )  =\Pi (
k_{-},k_{+},m,1110 )  ,\nonumber\\
\nu^{\prime} (  p_{4m+3,3} )      =\Pi (  k_{+},k_{-}
,m,1101 )  ,\qquad \nu^{\prime} (  p_{4m+3,4} )  =\Pi (
k_{+},k_{-},m,1110 )  .\label{nrm3}
\end{gather}
These formulae, together with the harmonic decomposition of polynomials, prove
that $ \langle \cdot,\cdot \rangle _{\tau}$ is positive-def\/inite for
$-\frac{1}{2}<k_{+},k_{-}<\frac{1}{2}$ (equivalently $-\frac{1}{2}<\pm
k_{0}\pm k_{1}<\frac{1}{2}$) (this is a special case of \cite[Proposition~4.4]{Etingof/Stoica:2009}).

These norm results are used to analyze the representation of the rational
Cherednik algeb\-ra~$\mathcal{A}_{\kappa}$ on $\mathcal{P}_{V}$ for arbitrary
parameter values. For f\/ixed $k_{0}$, $k_{1}$ the radical in~$\mathcal{P}_{V}$ is
the subspace
\[
\mathrm{rad}_{V} (  k_{0},k_{1} )  := \{  p: \langle
p,q \rangle _{\tau}=0~\forall\, q\in\mathcal{P}_{V} \}  .
\]
The radical is an $\mathcal{A}_{\kappa}$-module and the representation is
called \textit{unitary} if the form $ \langle \cdot,\cdot \rangle
_{\tau}$ is positive-def\/inite on $\mathcal{P}_{V}/\mathrm{rad}_{V} (
k_{0},k_{1} )$. By Proposition~\ref{Hbasis} $ \{   \vert
x \vert ^{2m}p_{n,i} \}  $ (with $m,n\in\mathbb{N}_{0}$, $1\leq
i\leq4$, except $1\leq i\leq2$ when $n=0$) is a basis for~$\mathcal{P}_{V}$
for any parameter values (see formula~(\ref{Hdecomp})).

\begin{proposition}
For fixed $k_{0}$, $k_{1}$ the set $ \{   \vert x \vert ^{2m}p_{n,i}:\nu (  p_{n,i} )  =0 \}  $ is a basis for $\mathrm{rad}
_{V} (  k_{0},k_{1} )$.
\end{proposition}

\begin{proof}
Suppose $p=\!\!\sum\limits_{m,n,i}\!\!a_{m,n,i} \vert x \vert ^{2m}p_{n,i}%
\!\in\!\mathrm{rad}_{V} (  k_{0},k_{1} )  $ then $ \langle
p, \vert x \vert ^{2\ell}p_{k,j} \rangle _{\tau}\!=\!a_{\ell
,k,j} 
4^{\ell}\ell! (  k{+}1 )  _{\ell}\nu (  p_{k,j} )  $
by~(\ref{xsqform}), thus $a_{\ell,k,j}\neq0$ implies $\nu (  p_{k,j} )
=0$.
\end{proof}

For a given value $k_{+}=\frac{1}{2}+m_{+}$ or $k_{-}=\frac{1}{2}+m_{-}$ (or
both) with $m_{+},m_{-}\in\mathbb{Z}$ the polynomials with $\nu\left(
p_{k,j}\right)  =0$ can be determined from the above norm formulae. For
example let $k_{+}=-\frac{1}{2}$ then $\mathrm{rad}_{V} \big(  k_{0},-\frac
{1}{2}-k_{0} \big)  $ contains $p_{n,i}$ for $ (  n,i )  $ in the
list $(1,2 )  $, $(  n,1 )$, $(  n,2)$
 ($n\geq2$ and $n\equiv0,1\operatorname{mod}4$), $ (  n,3 )$, $(
n,4 )  $ ($n\geq2$ and $n\equiv2,3\operatorname{mod}4$). In
particular each point on the boundary of the square $-\frac{1}{2}<k_{0}\pm
k_{1}<\frac{1}{2}$ corresponds to a unitary representation of~$\mathcal{A}_{\kappa}$.

\section{The reproducing kernel}\label{section4}

In the $\tau=1$ setting there is a (``\textit{Dunkl}'') kernel $E(x,y)  $ which
satisf\/ies $\mathcal{D}_{i}^{(x)  }E(x,y)
=y_{i}E(x,y)  $, $E(x,y)  =E (  y,x )  $
and $ \langle E (  \cdot,y )  ,p \rangle _{1}=p (
y )  $ for any polynomial~$p$. We show there exists such a~function in
this $B_{2}$-setting which is real-analytic in its arguments, provided
$\frac{1}{2}\pm k_{0}\pm k_{1}\notin\mathbb{Z}$. This kernel takes values in
$V\otimes V$; for notational convenience we will use expressions of the form
$\sum\limits_{i=1}^{2}\sum\limits_{j=1}^{2}f_{ij}(x,y)  s_{i}t_{j}$ where each
$f_{ij}(x,y)  $ is a polynomial in $x_{1}$, $x_{2}$, $y_{1}$, $y_{2}$
(technically, we should write $s_{i}\otimes t_{j}$ for the basis elements in
$V\otimes V$). For a polynomial $f(x)  =f_{1}(x)
t_{1}+f_{2}(x)  t_{2}$ let~$f(y)  ^{\ast}$ denote
$f_{1}(y)  s_{1}+f_{2}(y)  s_{2}$. The kernel~$E$ is
def\/ined as a sum of terms like $\frac{1}{\nu (  p_{n,i} )  }
p_{n,i}(y)  ^{\ast}p_{n,i}(x)  $; by the
orthogonality relations (for $m,n=0,1,2,\ldots$ and $1\leq i,j\leq4$)
\[
\frac{1}{\nu (  p_{n,i} )  }p_{n,i}(y)  ^{\ast
} \langle p_{n,i},p_{m,j} \rangle _{\tau}=\delta_{mn}\delta
_{ij}p_{n,i}(y)  ^{\ast}.
\]
By formulae (\ref{deltaF}) and (\ref{xsqform})%
\[
 \big\langle  \vert x \vert ^{2a}p_{n,i}(x)  , \vert
x \vert ^{2b}p_{m,j}(x)   \big\rangle_{\tau}=\delta
_{nm}\delta_{ab}\delta_{ij}4^{a}a! (  n+1 )  _{a} \nu (
p_{n,i} )  .
\]
We will f\/ind upper bounds on $ \{  p_{n,i} \}  $ and lower bounds on
$ \{  \nu (  p_{n,i} )   \}  $ in order to establish
convergence properties of $E(x,y)  $. For $u\in\mathbb{R}$ set%
\[
d(u)  :=\min_{m\in\mathbb{Z}}\left\vert u+\frac{1}{2}
+m\right\vert .
\]
The condition that $\nu (  p_{n,i} )  \neq0$ for all $n\geq1,1\leq
i\leq4$ is equivalent to $d (  k_{+} )  d (  k_{-} )  >0$,
since each factor in the numerator of $\Pi$ is of the form $\frac{1}{4}%
+\frac{1}{2} (  \pm k_{0}\pm k_{1} )  +m$ or $\frac{3}{4}+\frac{1}%
{2} (  \pm k_{0}\pm k_{1} )  +m$ for $m=0,1,2,\ldots$.

\begin{definition}
Let $P_{0}(x,y)  :=s_{1}t_{1}+s_{2}t_{2}$ and for $n\geq1$ let%
\[
P_{n}(x,y)  :=\sum_{i=1}^{4}\frac{1}{\nu (  p_{n,i} )
}p_{n,i}(y)  ^{\ast}p_{n,i}(x)  .
\]
For $n\geq0$ let
\[
E_{n}(x,y)  :=\sum_{0\leq m\leq n/2}\frac{\left(  x_{1}^{2}%
+x_{2}^{2}\right)  ^{m}\left(  y_{1}^{2}+y_{2}^{2}\right)  ^{m}}%
{4^{m}m!\left(  n-2m+1\right)  _{m}}P_{n-2m}(x,y)  .
\]
\end{definition}

\begin{proposition}
If $k_{0}$, $k_{1}$ are generic or $d (  k_{+} )  d (  k_{-} )
>0$, and $n\geq1$ then $ \langle E_{n} (  \cdot,y )
,f \rangle _{\tau}=f(y)  ^{\ast}$ for each $f\in
\mathcal{P}_{V,n}$, and $\mathcal{D}_{i}^{(x)  }E_{n} (
x,y )  =y_{i}E_{n-1}(x,y)  $ for~$i=1,2$.
\end{proposition}

\begin{proof}
By hypothesis on $k_{0}$, $k_{1}$ the polynomial $E_{n}(x,y)  $
exists, and is a rational function of~$k_{0}$,~$k_{1}$. The reproducing property
is a consequence of the orthogonal decomposition $\mathcal{P}_{V,n}
=\sum\limits_{0\leq m\leq n/2}\oplus \vert x \vert ^{2m}\mathcal{H}
_{V,\kappa,n-2m}$. For the second part, suppose $f\in\mathcal{P}_{V,n-1}$,
then%
\begin{gather*}
\big\langle \mathcal{D}_{i}^{(x)  }E_{n}(x,y)
,f(x)  \big\rangle _{\tau}= \langle E_{n} (
x,y )  ,x_{i}f(x)   \rangle _{\tau}
=y_{i}f(y)  ^{\ast}=y_{i} \langle E_{n-1}(x,y)
,f(x)  \rangle _{\tau}.
\end{gather*}
If $-\frac{1}{2}<\pm k_{0}\pm k_{1}<\frac{1}{2}$ the bilinear form is
positive-def\/inite which implies $\mathcal{D}_{i}^{(x)  }
E_{n}(x,y)  =y_{i}E_{n-1}(x,y)  $. This is an
algebraic (rational) relation which holds on an open set of its arguments and
is thus valid for all $k_{0}$, $k_{1}$ except for the poles of~$E_{n}$.
\end{proof}

Write $E_{n}(x,y)  =\sum\limits_{i,j=1}^{2}E_{n,ij}(x,y)
s_{i}t_{j}$ then the equivariance relation becomes $E_{n,ij} (
xw,yw )  = (  w^{-1}E_{n,\cdot\cdot}(x,y)  w )
_{ij}$ for each $w\in W (  B_{2} )  $. Next we show that $\sum\limits
_{n=0}^{\infty}E_{n}(x,y)  $ converges to a real-entire function
in $x,y$. For $p(x,t)  =p_{1}(x)  t_{1}+p_{2} (
x )  t_{2}\in\mathcal{P}_{V}$ let $\beta (  p )   (
x )  =p_{1}(x)  ^{2}+p_{2}(x)  ^{2}$.

\begin{lemma}
Suppose $m\geq0$ and $i=1$ or $3$ then%
\begin{gather*}
   \beta (  p_{2m+1,i} )  (x)  +\beta (
p_{2m+1,i+1} )  (x)
   = (  x_{1}^{2}+x_{2}^{2} )   (  \beta (  p_{2m,i} )
(x)  +\beta (  p_{2m,i+1} )  (x)  ).
\end{gather*}
\end{lemma}

\begin{proof}
By Def\/inition~\ref{pdef} $p_{2m+1,i}=x_{1}p_{2m,i}+x_{2}p_{2m,i+1}$ and
$p_{2m+1,i+1}=-x_{2}p_{2m,i}+x_{1}p_{2m,i+1}$. The result follows from direct
calculation using these formulae.
\end{proof}

\begin{lemma}
Suppose $m\geq1$ then
\begin{gather*}
\beta (  p_{2m,1} )  (x)  +\beta (  p_{2m,2} )
(x)\\
{}=\big(  x_{1}^{2}+x_{2}^{2}\big) \left(  \left(  \frac{2m-1+2k_{-}}{2m-1}\right)  ^{2}\beta (
p_{2m-1,3} )  (x)  +\left(  \frac{2m-1-2k_{-}}{2m-1}\right)
^{2}\beta (  p_{2m-1,4} )  (x)  \right)  ,\\
\beta (  p_{2m,3} )  (x)  +\beta (  p_{2m,4} )
(x)\\
=\big(  x_{1}^{2}+x_{2}^{2}\big)  \left(  \left(  \frac{2m-1+2k_{+}}{2m-1}\right)  ^{2}\beta (
p_{2m-1,1} )  (x)  +\left(  \frac{2m-1-2k_{+}}{2m-1}\right)
^{2}\beta (  p_{2m-1,2} )  (x)  \right)  .
\end{gather*}
\end{lemma}

\begin{proof}
This follows similarly as the previous argument.
\end{proof}

\begin{proposition}
Suppose $m\geq0$ and $i=1$ or $3$ then%
\begin{gather*}
\beta (  p_{2m,i} )  (x)  +\beta (  p_{2m,i+1}
 )  (x)      \leq2\left(  \frac{ (  1/2+ \vert
k_{0} \vert + \vert k_{1} \vert  )  _{m}}{ (
1/2 )  _{m}} \right)  ^{2}\big(  x_{1}^{2}+x_{2}^{2}\big)  ^{2m},\\
\beta (  p_{2m+1,i} )  (x)  +\beta (  p_{2m+1,i+1}  )  (x)     \leq2\left(  \frac{ (  1/2+ \vert
k_{0} \vert + \vert k_{1} \vert  )  _{m}}{ (
1/2 )  _{m}} \right)  ^{2}\big(  x_{1}^{2}+x_{2}^{2}\big)  ^{2m+1}.
\end{gather*}
\end{proposition}

\begin{proof}
Use induction, the lemmas, and the inequality $\big(  \frac{2m-1\pm2k_{0}%
\pm2k_{1}}{2m-1}\big)  ^{2}\leq \big(  \frac{m-1/2+ \vert k_{0}%
 \vert + \vert k_{1} \vert }{m-1/2}\big)  ^{2}$. The beginning
step is $\beta (  p_{1,i} )  (x)  +\beta (
p_{1,i+1} )  (x)  =2\big(  x_{1}^{2}+x_{2}^{2}\big)  $.
\end{proof}

By Stirling's formula
\begin{gather*}
\frac{(  1/2+ \vert k_{0} \vert + \vert k_{1} \vert
 )  _{m}}{ (  1/2 )  _{m}}     =\frac{\Gamma (  1/2 )
}{\Gamma (  1/2+ \vert k_{0} \vert + \vert k_{1} \vert
 )  }\frac{\Gamma (  1/2+ \vert k_{0} \vert + \vert
k_{1} \vert +m )  }{\Gamma (  1/2+m )  }\\
\hphantom{\frac{(  1/2+ \vert k_{0} \vert + \vert k_{1} \vert
 )  _{m}}{ (  1/2 )  _{m}}  }{}
  \sim\frac{\Gamma (  1/2 )  }{\Gamma (  1/2+ \vert
k_{0} \vert + \vert k_{1} \vert  )  }m^{ \vert
k_{0} \vert + \vert k_{1} \vert }
\end{gather*}
as $m\rightarrow\infty$. For the purpose of analyzing lower bounds on
$ \vert \nu^{\prime} (  p_{n,i} )   \vert $ we consider an
inf\/inite product.

\begin{lemma}
Suppose $u>0$ then the infinite product
\[
\omega (  u;z )
:=\prod\limits_{n=0}^{\infty}\left(  1-\frac{z^{2}}{(  u+n )  ^{2}
}\right)
\]
converges to an entire function of $z\in\mathbb{C}$, satisfies
\[
\omega (  u;z )  =\frac{\Gamma(u)  ^{2}}{\Gamma (
u+z )  \Gamma (  u-z )  },
\]
and has simple zeros at $\pm (  u+n )$, $n=0,1,2,\ldots$.
\end{lemma}

\begin{proof}
The product converges to an entire function by the comparison test:
$\sum\limits_{n=1}^{\infty}\frac{ \vert z \vert ^{2}}{ (  u+n )
^{2}}<\infty$ (this means that the partial products converge to a nonzero
limit, unless one of the factors vanishes). Suppose that $ \vert
z \vert <u$ then $\operatorname{Re} (  u\pm z )  \geq
u- \vert z \vert >0$ and
\begin{gather*}
\prod\limits_{n=0}^{\infty}\left(  1-\frac{z^{2}}{ (  u+n )  ^{2}}\right)
 =\lim_{m\rightarrow\infty}\frac{ (  u+z )  _{m} (
u-z )  _{m}}{(u)  _{m}^{2}}\\
\qquad{}
=\frac{\Gamma(u)  ^{2}}{\Gamma (  u+z )  \Gamma (
u-z )  }\lim_{m\rightarrow\infty}\frac{\Gamma (  u+z+m )
\Gamma (  u-z+m )  }{\Gamma (  u+m )  ^{2}}
=\frac{\Gamma(u)  ^{2}}{\Gamma (  u+z )  \Gamma (
u-z )  },
\end{gather*}
by Stirling's formula. The entire function $\omega (  u;\cdot )  $
agrees with the latter expression (in $\Gamma$) on an open set in $\mathbb{C}$, hence for all $z$.
\end{proof}

\begin{corollary}\qquad
\begin{enumerate}\itemsep=0pt
\item[$1.$] If $n\equiv0,1\operatorname{mod}4$ and $i=1,2$ or $n\equiv
2,3\operatorname{mod}4$ and $i=3,4$ then%
\[
\lim_{n\rightarrow\infty}\nu^{\prime} (  p_{n,i} )  =\omega \left(
\frac{1}{4};\frac{k_{+}}{2}\right)  \omega\left(  \frac{3}{4};\frac{k_{-}}%
{2}\right)  .
\]
\item[$2.$] If $n\equiv0,1\operatorname{mod}4$ and $i=3,4$ or $n\equiv
2,3\operatorname{mod}4$ and $i=1,2$ then
\[
\lim_{n\rightarrow\infty} \nu^{\prime} (  p_{n,i} )  =\omega\left(
\frac{1}{4};\frac{k_{-}}{2}\right)  \omega\left(  \frac{3}{4};\frac{k_{+}}%
{2}\right)  .
\]
\end{enumerate}
\end{corollary}

\begin{proof}
The statements follow from formulae (\ref{norms1}), (\ref{nrm2}) and
(\ref{nrm3}).
\end{proof}

For $a\in\mathbb{R}$ with $a+\frac{1}{2}\notin\mathbb{Z}$ let $\left(  \left(
a\right)  \right)  $ denote the nearest integer to $a$.

\begin{lemma}
For each $ (  n_{1},n_{2} )  \in\mathbb{Z}^{2}$ there is a constant
$C (  n_{1},n_{2} )  $ such that $ (   (  k_{+} )
 )  =n_{1}$ and $ (   (  k_{-} )   )  =n_{2}$ implies
$ \vert \nu^{\prime} (  p_{n,i} )   \vert \geq C (
n_{1},n_{2} )  d (  k_{+} )  d (  k_{-} )  $.
\end{lemma}

\begin{proof}
It suf\/f\/ices to consider Case~1  in the corollary (interchange $k_{+}$ and
$k_{-}$ to get Case~2). The limit function has zeros at $k_{+}=\pm\frac
{4m+1}{2}$ and $k_{-}=\pm\frac{4m+3}{2}$ for $m=0,1,2,\ldots$. For each
$n_{1}\neq0$ and $n_{2}\neq0$ there are unique nearest zeros $k_{+}=z_{1}$ and
$k_{-}=z_{2}$ respectively; for example if $n_{1}$ is odd and $n_{1}\geq1$
then $z_{1}=n_{1}-\frac{1}{2}$; and if $n_{2}$ is even and $n_{2}\leq-2$ then
$z_{2}=n_{2}+\frac{1}{2}$. Consider the entire function
\[
f\left(  k_{+},k_{-}\right)  =\frac{1}{ (  z_{1}-k_{+} )   (
z_{2}-k_{-} )  }\omega\left(  \frac{1}{4};\frac{k_{+}}{2}\right)
\omega\left(  \frac{3}{4};\frac{k_{-}}{2}\right)  .
\]
If $n_{1}=0$ then replace the factor $ (  z_{1}-k_{+} )  $ by
$\big(  \frac{1}{4}-k_{+}^{2}\big)  $, and if $n_{2}=0$ replace $\left(
z_{2}-k_{-}\right)  $ by~$1$. In each of these cases the quotient is an entire
function with no zeros in $ \vert n_{1}-k_{+} \vert \leq\frac{1}{2}$
and $ \vert n_{2}-k_{-} \vert \leq\frac{1}{2}$. Thus there is a lower
bound $C_{1}$ in absolute value for all the partial products of~$f$, valid
for all $k_{+}$, $k_{-}$ in this region. The expressions for $\nu^{\prime} (
p_{n,i} )  $ (see formulae~(\ref{norms1}),~(\ref{nrm2}) and~(\ref{nrm3}))
involve terms in $ \{  \varepsilon_{i} \}  $ but these do not af\/fect
the convergence properties (and note $0\leq\sum\limits_{i=1}^{4}\varepsilon_{i}\leq3$
in each case). Since $ \vert z_{1}-k_{+} \vert \geq d (
k_{+} )  $ and $ \vert z_{2}-k_{-} \vert \geq d (
k_{-} )  $ we f\/ind that the partial products of $ (  z_{1}
-k_{+} )   (  z_{2}-k_{-} )  f (  k_{+},k_{-} )  $,
that is, the values of $\nu (  p_{n,i} )  / (  2^{n}n! )  $,
are bounded below by $d (  k_{+} )  d (  k_{-} )  C_{1}$ in
absolute value. In case $n_{1}=0$ the factor $\big(  \frac{1}{4}-k_{+}%
^{2}\big)  =d (  k_{+} )   (  1-d (  k_{+} )   )
\geq\frac{1}{2}d (  k_{+} )  $ (and in the more trivial case
$n_{2}=0$ one has $1>\frac{1}{2}\geq d (  k_{-} )  $).
\end{proof}

\begin{theorem}
For a fixed $k_{0},k_{1}\in\mathbb{R}$ satisfying $d (  k_{+} )
d (  k_{-} )  >0$ the series $E(x,y)  :=\sum\limits_{n=0}^{\infty}E_{n}(x,y)  $ converges absolutely and uniformly
on $\big\{  (x,y)  \in\mathbb{R}^{4}: \vert x \vert
^{2}+ \vert y \vert ^{2}\leq R^{2}\big\}  $ for any $R>0$.
\end{theorem}

\begin{proof}
We have shown there is a constant $C^{\prime}>0$ such that
\[
 \vert P_{n}(x,y)   \vert \leq C^{\prime}\frac
{ \vert x \vert ^{n} \vert y \vert ^{n}}{2^{n}n!}\left(
\frac{n}{2}\right)  ^{ \vert k_{0} \vert + \vert k_{1} \vert
},
\]
and thus the series
\begin{gather*}
E(x,y)      =\sum_{n=0}^{\infty}\sum_{0\leq m\leq n/2}
\frac{ \vert x \vert ^{2m} \vert y \vert ^{2m}}
{4^{m}m! (  n-2m+1 )  _{m}}P_{n-2m}(x,y)
   =\sum_{l=0}^{\infty}P_{l}(x,y)  \sum_{m=0}^{\infty}
\frac{ \vert x \vert ^{2m} \vert y \vert ^{2m}}
{4^{m}m! (  l+1 )  _{m}}
\end{gather*}
converges uniformly for $ \vert x \vert ^{2}+ \vert y \vert
^{2}<R^{2}$.
\end{proof}

\begin{corollary}
Suppose $f\in\mathcal{P}_{V}$ then $ \langle E (  \cdot,y )
,f \rangle _{\tau}=f(y)  ^{\ast}$; also $\mathcal{D}_{i}^{(x)  }E(x,y)  =y_{i}E(x,y)  $ for
$i=1,2$.
\end{corollary}

As in the scalar ($\tau=1$) theory the function $E(x,y)  $ can be
used to def\/ine a generalized Fourier transform.

\section{The Gaussian-type weight function}\label{section5}

In this section we use vector notation for $\mathcal{P}_{V}$: $f (
x )  = (  f_{1}(x)  ,f_{2}(x)   )  $
for the previous $f_{1}(x)  t_{1}+f_{2}(x)  t_{2}$.
The action of $W$ is written as $ (  wf )  (x)
=f (  xw )  w^{-1}$. We propose to construct a $2\times2$
positive-def\/inite matrix function $K(x)  $ on $\mathbb{R}^{2}$
such that%
\[
\langle f,g\rangle _{G}=\ \int_{\mathbb{R}
^{2}}f(x)  K(x)  g(x)  ^{T}
e^{-\left\vert x\right\vert ^{2}/2}dx,\qquad \forall\, f,g\in\mathcal{P}_{V},
\]
and with the restriction $ \vert k_{0}\pm k_{1} \vert <\frac{1}{2}$;
the need for this was demonstrated in the previous section. The two necessary
algebraic conditions are (for all $f,g\in\mathcal{P}_{V}$)
\begin{gather}
 \langle f,g \rangle _{G}     = \langle wf,wg \rangle
_{G}, \qquad w\in W,\label{relGw}\\
 \langle  (  x_{i}-\mathcal{D}_{i} )  f,g \rangle _{G}
= \langle f,\mathcal{D}_{i}g \rangle _{G}, \qquad i=1,2. \label{relGDx}
\end{gather}
We will assume $K$ is dif\/ferentiable on $\Omega:= \{  x\in\mathbb{R} ^{2}:x_{1}x_{2} (  x_{1}^{2}-x_{2}^{2} )  \neq0 \}  $. The
integral formula is def\/ined if $K$ is integrable, but for the purpose of
dealing with the singularities implicit in $\mathcal{D}_{i}$ we introduce the
region%
\[
\Omega_{\varepsilon}:= \{  x:\min (   \vert x_{1} \vert
, \vert x_{2} \vert , \vert x_{1}-x_{2} \vert , \vert
x_{1}+x_{2} \vert  )  \geq\varepsilon \}
\]
for $\varepsilon>0$. The fundamental region of $W$ corresponding to $R_{+}$ is
$\mathcal{C}_{0}= \{  x:0<x_{2}<x_{1} \}  $.

Condition (\ref{relGw}) implies, for each $w\in W$, $f,g\in\mathcal{P}_{V}$
that
\begin{gather*}
\int_{\mathbb{R}
^{2}}f (  xw )  w^{-1}K(x)  wg (  xw )
^{T}e^{- \vert x \vert ^{2}/2}dx     =\int_{\mathbb{R}
^{2}}f(x)  w^{-1}K \big(  xw^{-1} \big)  wg(x)
^{T}e^{- \vert x \vert ^{2}/2}dx\\
\hphantom{\int_{\mathbb{R}
^{2}}f (  xw )  w^{-1}K(x)  wg (  xw )
^{T}e^{- \vert x \vert ^{2}/2}dx }{}
   =\int_{\mathbb{R}
^{2}}f(x)  K(x)  g(x)  ^{T}
e^{- \vert x \vert ^{2}/2}dx.
\end{gather*}
For the second step change the variable from $x$ to $xw^{-1}$ (note
$w^{T}=w^{-1}$). Thus we impose the condition
\[
K (  xw )  =w^{-1}K(x)  w.
\]
This implies that it suf\/f\/ices to determine $K$ on the fundamental region
$\mathcal{C}_{0}$ and then extend to all of $\Omega$ by using this formula.
Set $\partial_{i}:=\frac{\partial}{\partial x_{i}}$. Recall for the scalar
situation that the analogous weight function is $h_{\kappa}(x)
^{2}$ where $h_{\kappa}(x)  =\prod\limits_{v\in R_{+}} \vert
 \langle x,v \rangle  \vert ^{\kappa(v)  }$ and
satisf\/ies%
\[
\partial_{i}h_{\kappa}(x)  =\sum_{v\in R_{+}}\kappa (
v )  \frac{v_{i}}{ \langle x,v \rangle }h_{\kappa} (
x )  ,\qquad i=1,2.
\]
Start by solving the equation
\begin{gather}
\partial_{i}L(x)  =\sum_{v\in R_{+}}\kappa(v)
\frac{v_{i}}{ \langle x,v \rangle }L(x)  \sigma
_{v}, \qquad i=1,2, \label{eqnL}
\end{gather}
for a $2\times2$ matrix function $L$ on $\mathcal{C}_{0}$, extended by
\[
L (  xw )  =L(x)  w
\]
(from the facts that $w^{-1}\sigma_{v}w=\sigma_{vw}$ and $\kappa (
vw )  =\kappa(v)  $ it follows that equation~(\ref{eqnL}) is
satisf\/ied on all of~$\Omega$) and set
\[
K(x)  :=L(x)  ^{T}L(x)  ,
\]
with the result that~$K$ is positive-semidef\/inite and (note $\sigma_{v}^{T}=\sigma_{v}$)
\begin{gather}
\partial_{i}K(x)  =\sum_{v\in R_{+}}\kappa(v)
\frac{v_{i}}{ \langle x,v \rangle } (  \sigma_{v}K (
x )  +K(x)  \sigma_{v} )  . \label{eqnK}
\end{gather}
Then for $f,g\in\mathcal{P}_{V}$ (and $i=1,2$) we f\/ind
\begin{gather*}
\begin{split}
& -\partial_{i}\big(  f(x)  K(x)  g(x)
^{T}e^{- \vert x \vert ^{2}/2}\big)  e^{ \vert x \vert
^{2}/2}=
x_{i}f(x)  K(x)  g(x)  ^{T}- (
\partial_{i}f(x)   )  K(x)  g(x)
^{T}\\
& \qquad{}-f(x)  ^{T}K(x)   (  \partial_{i}g (
x )   )  ^{T}
-\sum_{v\in R_{+}}\kappa(v)  \frac{v_{i}}{ \langle
x,v \rangle }\big\{  f(x)  \sigma_{v}K(x)
g(x)  ^{T}+f(x)  K(x)  \sigma
_{v}g(x)  ^{T}\big\}.
\end{split}
\end{gather*}
Consider the second necessary condition (\ref{relGDx}) $ \langle  (
x_{i}-\mathcal{D}_{i} )  f,g \rangle _{G}- \langle
f,\mathcal{D}_{i}g \rangle _{G}=0$, that is, the following integral must
vanish
\begin{gather*}
   \int_{\Omega_{\varepsilon}}\big\{   (   (  x_{i}-\mathcal{D}%
_{i} )  f(x)   )  K(x)  g(x)
^{T}-f(x)  K(x)   (  \mathcal{D}_{i}g (
x\ )   )  ^{T}\big\}  e^{- \vert x \vert ^{2}/2}dx\\
 \qquad{}  =-\int_{\Omega_{\varepsilon}}\partial_{i}\big(  f(x)
K(x)  g(x)  ^{T}e^{- \vert x \vert ^{2}
/2}\big)  dx\\
\qquad\quad{}  +\sum_{v\in R_{+}}\kappa(v)  \int_{\Omega_{\varepsilon}}
\frac{v_{i}}{ \langle x,v \rangle }\big\{  f (  x\sigma
_{v} )  \sigma_{v}K(x)  g(x)  ^{T}+f (
x )  K(x)  \sigma_{v}g (  x\sigma_{v} )
^{T}\big\}  dx.
\end{gather*}
In the second part, for each $v\in R_{+}$ change the variable to $x\sigma_{v}$, then the numerator is invariant, because $\sigma_{v}K (  x\sigma
_{v} )  =K(x)  \sigma_{v}$, and $ \langle x\sigma
_{v},v \rangle =- \langle x,v \rangle $, and thus each term
vanishes (note $\Omega_{\varepsilon}$ is $W$-invariant). So establishing the
validity of the inner product formula reduces to showing $\lim
\limits_{\varepsilon\rightarrow0_{+}}\int_{\Omega_{\varepsilon}}\partial
_{i}\big(  f(x)  K(x)  g(x)
^{T}e^{-\ \vert x \vert ^{2}/2}\big)  dx=0$ for $i=1,2$. By the
polar identity it suf\/f\/ices to prove this for $g=f$. Set $Q(x)
=f(x)  K(x)  f(x)  ^{T}e^{- \vert
x \vert ^{2}/2}$.

\begin{figure}[t]
\centering
\includegraphics[width=3.0in]{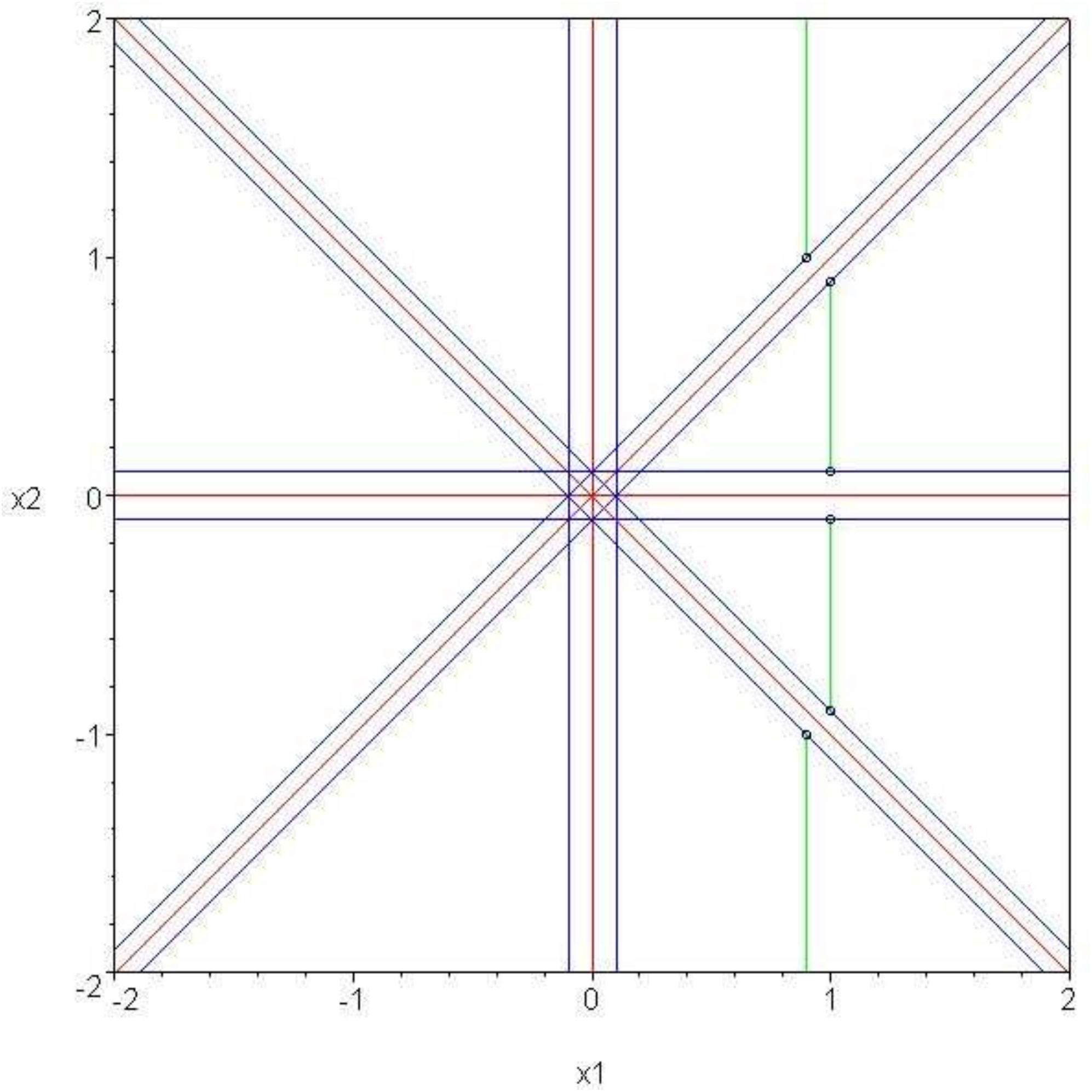}
\caption{Region of integration.}\label{omega}
\end{figure}

By symmetry ($\sigma_{12}^{+}\mathcal{D}_{1}\mathcal{\sigma}_{12}
^{+}=\mathcal{D}_{2}$) it suf\/f\/ices to prove the formula for $i=2$. Consider
the part of~$\Omega_{\varepsilon}$ in $ \{  x_{1}>0 \}  $ as the
union of $ \{  x:\varepsilon< \vert x_{2} \vert <x_{1}
-\varepsilon \}  $ and $ \{  x:\varepsilon<x_{1}< \vert
x_{2} \vert -\varepsilon \}  $ (with vertices $ (
2\varepsilon,\pm\varepsilon )  $, $ (  \varepsilon,\pm2\varepsilon
 ) $ respectively). In the iterated integral evaluate the inner integral
over $x_{2}$ on the segments $ \{   (  x_{1},x_{2} )
:\varepsilon\leq \vert x_{2} \vert \leq x_{1}-\varepsilon \}  $
and $ \{   (  x_{1}-\varepsilon,x_{2} )  :x_{1}\leq \vert
x_{2} \vert  \}  $ for a f\/ixed $x_{1}>2\varepsilon$; obtaining (due
to the exponential decay)
\begin{gather*}
  - (  Q (  x_{1}-\varepsilon,x_{1} )  -Q (  x_{1}
,x_{1}-\varepsilon )   )  - (  Q (  x_{1},\varepsilon
 )  -Q (  x_{1},-\varepsilon )   ) \\
\qquad{}  - (  Q (  x_{1},-x_{1}+\varepsilon )  -Q (  x_{1}
-\varepsilon,-x_{1} )   )  .
\end{gather*}
See Fig.~\ref{omega} for a diagram of $\Omega_{\varepsilon}$ and a typical
inner integral.

By dif\/ferentiability
\begin{gather*}
h_{ij} (  x_{1},\varepsilon )  -h_{ij} (  x_{1},-\varepsilon
 )      =C_{1} (  x_{1} )  \varepsilon+O\big(  \varepsilon
^{2}\big)  ,\\
h_{ij} (  x_{1}-\varepsilon,x_{1} )  -h_{ij} (  x_{1} ,x_{1}-\varepsilon )      =C_{2} (  x_{1} )  \varepsilon
+O\big(  \varepsilon^{2}\big)  ,
\end{gather*}
where $h_{ij}(x)  =f_{i}(x)  f_{j} (
x )  e^{- \vert x \vert ^{2}/2}$ for $1\leq i,j\leq2$; the
factors $C_{1}$, $C_{2}$ depend on $x_{1}$ but there is a global bound
$ \vert C_{i} (  x_{1} )   \vert <C_{0}$ depending only on
$f$ because of the exponential decay. Thus the behavior of $K (
x_{1},\varepsilon )  $ and $K (  x_{1},x_{1}-\varepsilon )  $ is
crucial in analyzing the limit as $\varepsilon\rightarrow0_{+}$, for
$x_{2}=-x_{1},0,x_{1}$. It suf\/f\/ices to consider the fundamental region
$0<x_{2}<x_{1}$. At the edge $ \{  x_{2}=0 \}  $ corresponding to
$\sigma_{2}$
\[
K (  x_{1},-\varepsilon )  =K (   (  x_{1},\varepsilon
 )  \sigma_{2} )  =\sigma_{2}K (  x_{1},\varepsilon )
\sigma_{2},
\]
and
\begin{gather*}
 (  Q (  x_{1},\varepsilon )  -Q (  x_{1},-\varepsilon
 )   )  = (  h_{11} (  x_{1},\varepsilon )
-h_{11} (  x_{1},-\varepsilon )  )  K (  x_{1}
,\varepsilon )  _{11}\\
\qquad{} + (  h_{22} (  x_{1},\varepsilon )  -h_{22} (  x_{1}
,-\varepsilon )   )  K (  x_{1},\varepsilon )
_{22}-2 (  h_{12} (  x_{1},\varepsilon )  +h_{12} (
x_{1},-\varepsilon )   )  K (  x_{1},\varepsilon )  _{12}.
\end{gather*}
At the edge $ \{  x_{1}-x_{2}=0 \}  $ corresponding to $\sigma^+_{12}$
\[
K (  x_{1}-\varepsilon,x_{1} )  =K (   (  x_{1}%
,x_{1}-\varepsilon )  \sigma^+_{12} )  =\sigma^+_{12}K (
x_{1},x_{1}-\varepsilon )  \sigma^+_{12}.
\]
For conciseness set $x^{ (  0 )  }= (  x_{1},x_{1}-\varepsilon
 )  $ and $x^{ (  1 )  }= (  x_{1}-\varepsilon,x_{1} )
$. Then
\begin{gather*}
Q\big(  x^{(0)  }\big)  -Q\big(  x^{(1)
}\big)  =2\big(  h_{12}\big(  x^{(0)  }\big)
-h_{12}\big(  x^{(1)  }\big)  \big)  K\big(  x^{(0)  }\big)  _{12}\\
\qquad{} +\frac{1}{2}\big(  h_{11}\big(  x^{(0)  }\big)
-h_{11}\big(  x^{(1)  }\big)  +h_{22}\big(  x^{(0)  }\big)  -h_{22}\big(  x^{(1)  }\big)  \big)
\big(  K\big(  x^{(0)  }\big)  _{11}+K\big(  x^{(0)  }\big)  _{22}\big) \\
\qquad{} +\frac{1}{2}\big(  h_{11}\big(  x^{(0)  }\big)
+h_{11}\big(  x^{(1)  }\big)  -h_{22}\big(  x^{(0)  }\big)  -h_{22}\big(  x^{(1)  }\big)  \big)
\big(  K\big(  x^{(0)  }\big)  _{11}-K\big(  x^{(0)  }\big)  _{22}\big)  .
\end{gather*}
To get zero limits as $\varepsilon\rightarrow0_{+}$ some parts rely on the
uniform continuity of $h_{ij}$ and bounds on certain entries of~$K$, and the
other parts require
\begin{gather*}
\lim\limits_{\varepsilon\rightarrow0_{+}}K (  x_{1},\varepsilon )
_{12}    =0,\qquad
\lim\limits_{\varepsilon\rightarrow0_{+}} (  K (  x_{1},x_{1}
-\varepsilon )  _{11}-K (  x_{1},x_{1}-\varepsilon )
_{22} )      =0.
\end{gather*}
This imposes
various conditions on $K$ near the edges, as described above.

We turn to the solution of the system (\ref{eqnL}) and rewrite
\begin{gather*}
\partial_{1}L(x)      =L(x)  \left\{  \frac{k_{1}%
}{x_{1}}%
\begin{bmatrix}
-1 & 0\\
0 & 1
\end{bmatrix}
+\frac{2k_{0}x_{2}}{x_{1}^{2}-x_{2}^{2}}%
\begin{bmatrix}
0 & 1\\
1 & 0
\end{bmatrix}
\right\}  ,
\\
\partial_{2}L(x)      =L(x)  \left\{  \frac{k_{1}%
}{x_{2}}%
\begin{bmatrix}
1 & 0\\
0 & -1
\end{bmatrix}
-\frac{2k_{0}x_{1}}{x_{1}^{2}-x_{2}^{2}}%
\begin{bmatrix}
0 & 1\\
1 & 0
\end{bmatrix}
\right\}  .\nonumber
\end{gather*}
The ref\/lections $\sigma_{12}^{-}$ and $\sigma_{12}^{+}$ were combined into one
term: $\frac{ (  1,-1 )  }{x_{1}-x_{2}}-\frac{(1,1)}{x_{1}+x_{2}}
=\frac{ (  2x_{2},-2x_{1} )  }{x_{1}^{2}-x_{2}^{2}}$. Since $ (
x_{1}\partial_{1}+x_{2}\partial_{2} )  L(x)  =0$ we see that
$L$ is positively homogeneous of degree~$0$. Because of the homogeneity the
system can be transformed to an ordinary dif\/ferential system by setting
$u=\frac{x_{2}}{x_{1}}$. Then the system is transformed to
\[
\frac{d}{du}L(u)  =L(u)  \left\{  \frac{k_{1}}{u}
\begin{bmatrix}
1 & 0\\
0 & -1
\end{bmatrix}
-\frac{2k_{0}}{1-u^{2}}%
\begin{bmatrix}
0 & 1\\
1 & 0
\end{bmatrix}
\right\}  .
\]
It follows from this equation that $\frac{d}{du}\det L(u)  =0$.
Since the goal is to f\/ind a positive-def\/inite matrix $K$ we look for a
fundamental solution for~$L$, that is  $\det L(x)  \neq0$. If
$L(x)  $ is a solution then so is~$ML(x)  $ for any
nonsingular constant matrix. Thus $K(x)  =L(x)
^{T}M^{T}ML(x)  $ satisf\/ies the dif\/ferential equation~(\ref{eqnK}) and some other condition must be imposed to obtain the desired
(unique) solution for the weight function. The process starts by solving for a
row of~$L$, say $ (  f_{1}(u)  ,f_{2}(u)
 )  $, that is
\begin{gather*}
\frac{d}{du}f_{1}(u)      =\frac{k_{1}}{u}f_{1}(u)
-\frac{2k_{0}}{1-u^{2}}f_{2}(u)  ,\qquad
\frac{d}{du}f_{2}(u)      =-\frac{2k_{0}}{1-u^{2}}f_{1} (
u )  -\frac{k_{1}}{u}f_{2}(u)  .
\end{gather*}
A form of solutions can be obtained by computer algebra, then a desirable
solution can be verif\/ied. Set
\begin{gather*}
f_{1}(u)      = \vert u \vert ^{k_{1}} \big(
1-u^{2} \big)  ^{-k_{0}}g_{1} \big(  u^{2}\big)  ,\qquad
f_{2}(u)      = \vert u \vert ^{k_{1}}\big(
1-u^{2}\big)  ^{-k_{0}}ug_{2}\big(  u^{2}\big)  ,
\end{gather*}
then the equations become (with $s:=u^{2}$)
\begin{gather*}
 (  s-1 )  \frac{d}{ds}g_{1} (  s )      =k_{0}g_{1} (
s )  +k_{0}g_{2} (  s )  ,\\
s (  s-1 )  \frac{d}{ds}g_{2} (  s )      =k_{0}g_{1} (
s )  +\left\{  k_{0}s-\left(  \frac{1}{2}+k_{1}\right)   (
s-1 )  \right\}  g_{2} (  s )  ,
\end{gather*}
and the solution regular at $s=0$ is
\begin{gather*}
g_{1} (  s )      =F\left(  -k_{0},k_{1}+\frac{1}{2}-k_{0}
;k_{1}+\frac{1}{2};s\right)  ,\\
g_{2} (  s )      =-\frac{k_{0}}{k_{1}+\frac{1}{2}}F\left(
1-k_{0},k_{1}+\frac{1}{2}-k_{0};k_{1}+\frac{3}{2};s\right)  .
\end{gather*}
(We use $F$ to denote the hypergeometric function ${}_{2}F_{1}$; it is the only
type appearing here.) The verif\/ication uses two hypergeometric identities
(arbitrary parameters $a$, $b$, $c$ with $-c\notin\mathbb{N}_{0}$)
\begin{gather*}
 (  s-1 )  \frac{d}{ds}F (  a,b;c;s )      =-aF (
a,b;c;s )  +\frac{a (  c-b )  }{c}F (  a+1,b;c+1;s )
,\\
s (  s-1 )  \frac{d}{ds}F (  a+1,b;c+1;s )      = (
c-bs )  F (  a+1,b;c+1;s )  -cF (  a,b;c;s )  .
\end{gather*}

To get the other solutions we use the symmetry of the system, replace $k_{1}$
by $-k_{1}$ and interchange $f_{1}$ and $f_{2}$. We have a fundamental
solution $L(u)  $ given by
\begin{gather*}
L(u)  _{11}     = \vert u \vert ^{k_{1}}\big(
1-u^{2}\big)  ^{-k_{0}}F\left(  -k_{0},\frac{1}{2}-k_{0}+k_{1};k_{1}%
+\frac{1}{2};u^{2}\right), \nonumber\\ 
L(u)  _{12}     =-\frac{k_{0}}{k_{1}+\frac{1}{2}} \vert
u \vert ^{k_{1}}\big(  1-u^{2}\big)  ^{-k_{0}}uF\left(  1-k_{0}%
,\frac{1}{2}-k_{0}+k_{1};k_{1}+\frac{3}{2};u^{2}\right), \nonumber\\
L(u)  _{21}     =-\frac{k_{0}}{\frac{1}{2}-k_{1}} \vert
u \vert ^{-k_{1}}\big(  1-u^{2}\big)  ^{-k_{0}}uF\left(  1-k_{0}%
,\frac{1}{2}-k_{0}-k_{1};\frac{3}{2}-k_{1};u^{2}\right), \nonumber\\
L(u)  _{22}     = \vert u \vert ^{-k_{1}}\big(
1-u^{2}\big)  ^{-k_{0}}F\left(  -k_{0},\frac{1}{2}-k_{0}-k_{1};\frac{1}%
{2}-k_{1};u^{2}\right)  .\nonumber
\end{gather*}
Observe that $\lim\limits_{u\rightarrow0_{+}}\det L(u)  =1$, thus
$\det L(u)  =1$ for all $u$. We can write $L$ in the form%
\[
L (  x_{1},x_{2} )  =
\begin{bmatrix}
 \vert x_{2} \vert ^{k_{1}}x_{1}^{-k_{1}} & 0\\
0 &  \vert x_{2} \vert ^{-k_{1}}x_{1}^{k_{1}}%
\end{bmatrix}
\left[
\begin{matrix}
c_{11}(x)  & \frac{x_{2}}{x_{1}}c_{12}(x) \\
\frac{x_{2}}{x_{1}}c_{21}(x)  & c_{22}(x)
\end{matrix}
\right]  ,
\]
where each $c_{ij}$ is even in $x_{2}$ and is real-analytic in $0< \vert
x_{2} \vert <x_{1}$. In fact $L(x)  $ is thus def\/ined on
$\mathcal{C}_{0}\cup\mathcal{C}_{0}\sigma_{2}$. It follows that $K= (
ML )  ^{T}ML$ is integrable near $ \{  x_{2}=0 \}  $ if
$ \vert k_{1} \vert <\frac{1}{2}$, and $\lim\limits_{\varepsilon
\rightarrow0_{+}}K (  x_{1},\varepsilon )  _{12}=0$ exactly when
$M^{T}M$ is diagonal. The standard identity  \cite[15.8.1]{Olver/Lozier/Boisvert/Clark:2010}
\begin{gather}
F (  a,b;c;u )  = (  1-u )  ^{c-a-b}F (
c-a,c-b;c;u )  \label{hgf1}
\end{gather}
shows that there is a hidden symmetry for $k_{0}$
\begin{equation}
\begin{split}
& \big(  1-u^{2}\big)  ^{-k_{0}}F\left(  -k_{0},\frac{1}{2}-k_{0}%
+k_{1};k_{1}+\frac{1}{2};u^{2}\right) \\
&  \qquad{}   =\big(  1-u^{2}\big)  ^{k_{0}}F\left(  k_{0},\frac{1}{2}+k_{0}%
+k_{1};k_{1}+\frac{1}{2};u^{2}\right)  ,
\end{split}\label{hgf2}
\end{equation}
and similar equations for the other entries of $L$. Consider
\begin{gather*}
Q (  x_{1},-\varepsilon )  -Q (  x_{1},\varepsilon )
=K (  x_{1},\varepsilon )  _{11} (  h_{11} (  x_{1}%
,-\varepsilon )  -h_{11} (  x_{1},\varepsilon )   ) \\
\qquad{} -2K (  x_{1},\varepsilon )  _{12} (  h_{12} (  x_{1}%
,-\varepsilon )  +h_{12} (  x_{1},\varepsilon )   )
+K (  x_{1},\varepsilon )  _{22} (  h_{22} (  x_{1}%
,-\varepsilon )  -h_{22} (  x_{1},\varepsilon )  )  .
\end{gather*}
With diagonal $M^{T}M$ we f\/ind (note $x_{1}>2\varepsilon$ in the region, so
$\frac{\varepsilon}{x_{1}}<\frac{1}{2}$)%
\begin{gather*}
K (  x_{1},\varepsilon )  _{11}     =O\big(  x_{1}^{-2k_{1}%
}\varepsilon^{2k_{1}}\big)  +O\big(  x_{1}^{2k_{1}-2}\varepsilon^{2-2k_{1}%
}\big)  ,\\
K (  x_{1},\varepsilon )  _{12}     =O\big(  x_{1}^{-1-2k_{1}%
}\varepsilon^{1+2k_{1}}\big)  +O\big(  x_{1}^{-1+2k_{1}}\varepsilon
^{1-2k_{1}}\big)  ,\\
K (  x_{1},\varepsilon )  _{22}     =O\big(  x_{1}^{2k_{1}%
}\varepsilon^{-2k_{1}}\big)  +O\big(  x_{1}^{-2-2k_{1}}\varepsilon
^{2+2k_{1}}\big)  .
\end{gather*}
By the exponential decay we can assume that the double integral is over the
box $\max (   \vert x_{1} \vert , \vert x_{2} \vert
 )  \leq R$ for some $R<\infty$. Note $\int_{2\varepsilon}^{R} \big(
\frac{\varepsilon}{x_{1}} \big)  ^{\alpha}dx_{1}=\frac{1}{1+\alpha} \big(
R^{1-\alpha}\varepsilon^{\alpha}-2^{1-\alpha}\varepsilon \big)  $.

These bounds (recall $-\frac{1}{2}<k_{0},k_{1}<\frac{1}{2}$) show that%
\[
\lim\limits_{\varepsilon\rightarrow0_{+}}\int_{2\varepsilon}^{R} (
Q (  x_{1},-\varepsilon )  -Q (  x_{1},\varepsilon )
 )  dx_{1}=0.
\]

Next we analyze the behavior of this solution in a neighborhood of $t=1$, that
is  the ray $ \{   (  x_{1},x_{1} )  :x_{1}>0 \}  $. The
following identity  \cite[15.10.21]{Olver/Lozier/Boisvert/Clark:2010}  is
used
\begin{gather*}
F (  a+d,a+c;c+d;u )  =\frac{\Gamma (  c+d )  \Gamma (
-2a )  }{\Gamma (  c-a )  \Gamma (  d-a )  }F (
a+d,c+a;1+2a;1-u ) \\
\hphantom{F (  a+d,a+c;c+d;u )  =}{}
+ (  1-u )  ^{-2a}\frac{\Gamma (  c+d )  \Gamma (
2a )  }{\Gamma (  c+a )  \Gamma (  d+a )  }F (
d-a,c-a;1-2a;1-u )  .
\end{gather*}
Let%
\begin{gather*}
\eta (  k_{0},k_{1} )      :=\frac{\Gamma\left(  \frac{1}{2}
+k_{1}\right)  \Gamma (  2k_{0} )  }{\Gamma\left(  \frac{1}{2}
+k_{0}+k_{1}\right)  \Gamma (  k_{0} )  }
   =\frac{2^{2k_{0}-1}}{\sqrt{\pi}}\frac{\Gamma\left(  \frac{1}{2}
+k_{1}\right)  \Gamma\left(  \frac{1}{2}+k_{0}\right)  }{\Gamma\left(
\frac{1}{2}+k_{0}+k_{1}\right)  };
\end{gather*}
the latter equation follows from $\Gamma (  2a )  /\Gamma (
a )  =2^{2a-1}\Gamma (  a+\frac{1}{2} )  /\sqrt{\pi}$ (the
duplication formula). We will need the identity
\[
\eta (  k_{0},k_{1} )  \eta (  -k_{0},-k_{1} )  +\eta (
k_{0},-k_{1} )  \eta (  -k_{0},k_{1} )  =\frac{1}{2},
\]
proved by use of $\Gamma \big(  \frac{1}{2}+a \big)  \Gamma \big(  \frac
{1}{2}-a \big)  =\frac{\pi}{\cos\pi a}$. Thus%
\begin{gather*}
L(u)  _{11}     = \vert u \vert ^{k_{1}}\left\{\eta (
k_{0},k_{1} )  \big(  1-u^{2}\big)  ^{-k_{0}}F\left(  -k_{0},\frac
{1}{2}-k_{0}+k_{1};1-2k_{0};1-u^{2}\right)\right. \\
\left.\hphantom{L(u)  _{11}     =}{}
 +\eta (  -k_{0},k_{1} )  \big(  1-u^{2}\big)  ^{k_{0}}F\left(
k_{0},\frac{1}{2}+k_{0}+k_{1};1+2k_{0};1-u^{2}\right)  \right\},
\end{gather*}
by use of identity (\ref{hgf2}), and also
\begin{gather*}
L(u)  _{12}     =u \vert u \vert ^{k_{1}}\left\{-\eta (
k_{0},k_{1} )  \big(  1-u^{2}\big)  ^{-k_{0}} \,{} _{2}F_{1}\left(
1-k_{0},\frac{1}{2}-k_{0}+k_{1};1-2k_{0};1-u^{2}\right)\right. \\
\left.\hphantom{L(u)  _{12}     =}{}
 +\eta (  -k_{0},k_{1} )  \big(  1-u^{2}\big)  ^{k_{0}} \,
{}_{2}F_{1}\left(  1+k_{0},\frac{1}{2}+k_{0}+k_{1};1+2k_{0};1-u^{2}\right)  \right\}.
\end{gather*}
Transform again using (\ref{hgf1}) to obtain
\begin{gather*}
   F\!\left( \! 1-k_{0},\frac{1}{2}-k_{0}+k_{1};1-2k_{0};1-u^{2}\!\right)
   = \vert u \vert ^{-2k_{1}-1}F\!\left(\!  -k_{0},\frac{1}{2}%
-k_{0}-k_{1};1-2k_{0};1-u^{2}\!\right)  ,
\end{gather*}
and so on. All the hypergeometric functions we use are of one form and it is
convenient to introduce
\[
H (  a,b;s )  :=F\left(  a,a+b+\frac{1}{2};2a+1;1-s\right)  .
\]

By using similar transformations as for $L_{11}$ and $L_{12}$ we f\/ind%
\begin{gather*}
L(u)      =\left[
\begin{matrix}
\eta (  -k_{0},k_{1} )  & \eta (  k_{0},k_{1} ) \\
\eta (  -k_{0},-k_{1} )  & -\eta (  k_{0},-k_{1} )
\end{matrix}
\right]  \left[
\begin{matrix}
\big(  1-u^{2}\big)  ^{k_{0}} & 0\\
0 & \big(  1-u^{2}\big)  ^{-k_{0}}
\end{matrix}
\right] \\
\hphantom{L(u)      =}{}
 \times\left[
\begin{matrix}
H (  k_{0},k_{1};u^{2} )  & H (  k_{0},-k_{1};u^{2} ) \\
H (  -k_{0},k_{1};u^{2} )  & -H (  -k_{0},-k_{1};u^{2} )
\end{matrix}
\right]  \left[
\begin{matrix}
 \vert u \vert ^{k_{1}} & 0\\
0 &  \vert u \vert ^{-k_{1}}
\end{matrix}
\right]  .
\end{gather*}
Let $\Gamma$ denote the f\/irst matrix in the above formula. By direct
calculation we f\/ind that a necessary condition for
\[
\lim\limits_{\varepsilon
\rightarrow0_{+}} (  K (  x_{1},x_{1}-\varepsilon )
_{11}-K (  x_{1},x_{1}-\varepsilon )  _{22} )  =0,
\]
 where
$K(u)  = (  M^{\prime}\Gamma^{-1}L(u)   )
^{T} (  M^{\prime}\Gamma^{-1}L(u)   )  $ is that
$M^{\prime T}M^{\prime}$ is diagonal. The proof that
\[
\int_{2\varepsilon}%
^{R} (  Q (  x_{1},x_{1}-\varepsilon )  -Q (  x_{1}%
-\varepsilon,x_{1} )   )  dx_{1}\rightarrow0
\]
 is similar to the
previous case; a typical term in $K(u)  _{11}-K(u)
_{22}$ is
\[
\big(  1-u^{2}\big)  ^{2k_{0}}\big(  H\big(  k_{0},k_{1};u^{2}\big)
^{2} \vert u \vert ^{2k_{1}}-H\big(  k_{0},-k_{1};u^{2}\big)
^{2} \vert u \vert ^{-2k_{1}}\big)
\]
which is $O \big(   (  1-u )  ^{1+2k_{0}} \big)  $, tending to zero
as $u\rightarrow1_{-}$ for $ \vert k_{0} \vert <\frac{1}{2}$).

It remains to combine the two conditions: $K(u)  = (
ML(u)   )  ^{T} (  ML(u)   )  $ and
the previous one to f\/ind a unique solution for $M$: $M^{T}M= (
\Gamma^{-1} )  ^{T}D\Gamma^{-1}$ and both $M^{T}M$ and $D$ are
positive-def\/inite diagonal. Indeed
\begin{gather*}
D     =c\left[
\begin{matrix}%
\eta (  -k_{0},-k_{1} )  \eta (  -k_{0},k_{1} )  & 0\\
0 & \eta (  k_{0},k_{1} )  \eta (  k_{0},-k_{1} )
\end{matrix}
\right]  ,\\
M^{T}M     =2c\left[
\begin{matrix}%
\eta (  -k_{0},-k_{1} )  \eta (  k_{0},-k_{1} )  & 0\\
0 & \eta (  k_{0},k_{1} )  \eta (  -k_{0},k_{1} )
\end{matrix}
\right]  ,
\end{gather*}
for some $c>0$. Thus the desired matrix weight (in the region $\mathcal{C}
_{0}$, $0<x_{2}<x_{1}$, $u=x_{2}/x_{1}$) is given by
\begin{gather*}
K(u)  _{11}     =d_{1}L(u)  _{11}^{2}+d_{2}L (
u )  _{21}^{2},\\
K(u)  _{12}     =K(u)  _{21}=d_{1}L(u)
_{11}L(u)  _{12}+d_{2}L(u)  _{21}L(u)
_{22},\\
K(u)  _{22}     =d_{1}L(u)  _{12}^{2}+d_{2}L (
u )  _{22}^{2},
\end{gather*}
where
\begin{gather*}
d_{1}     =c\frac{\Gamma\left(  \frac{1}{2}-k_{1}\right)  ^{2}}{\cos\pi
k_{0}\Gamma\left(  \frac{1}{2}+k_{0}-k_{1}\right)  \Gamma\left(  \frac{1}%
{2}-k_{0}-k_{1}\right)  },\\
d_{2}     =c\frac{\Gamma\left(  \frac{1}{2}+k_{1}\right)  ^{2}}{\cos\pi
k_{0}\Gamma\left(  \frac{1}{2}+k_{0}+k_{1}\right)  \Gamma\left(  \frac{1}%
{2}-k_{0}+k_{1}\right)  }.
\end{gather*}
Also $\det K=d_{1}d_{2}=c^{2}\big(  1-\tan^{2}\pi k_{0}\tan^{2}\pi
k_{1}\big)  $. The expressions for $K_{ij}$ can be rewritten somewhat by
using the transformations~(\ref{hgf2}). Observe that the conditions $-\frac
{1}{2}<\pm k_{0}\pm k_{1}<\frac{1}{2}$ are needed for $d_{1},d_{2}>0$. The
normalization constant is to be determined from the condition $\int
_{\mathbb{R}^{2}}K(x)  _{11}e^{- \vert x \vert ^{2}%
/2}dx=1$. By the homogeneity of $K$ this is equivalent to evaluating
$\int_{-\pi}^{\pi}K (  \cos\theta,\sin\theta )  _{11}d\theta$ (or
$\int_{0}^{\pi/4} (  K_{11}+K_{22} )   (  \cos\theta,\sin
\theta )  d\theta$) (the integral looks dif\/f\/icult because~$K_{11}$
involves squares of hypergeometric functions with argument~$\tan^{2}\theta$).
Numerical experiments suggest the following conjecture for the normalizing constant
\[
c=\frac{\cos\pi k_{0}\cos\pi k_{1}}{2\pi}.
\]

\begin{figure}[t]
\centering
\includegraphics[width=2.6in]{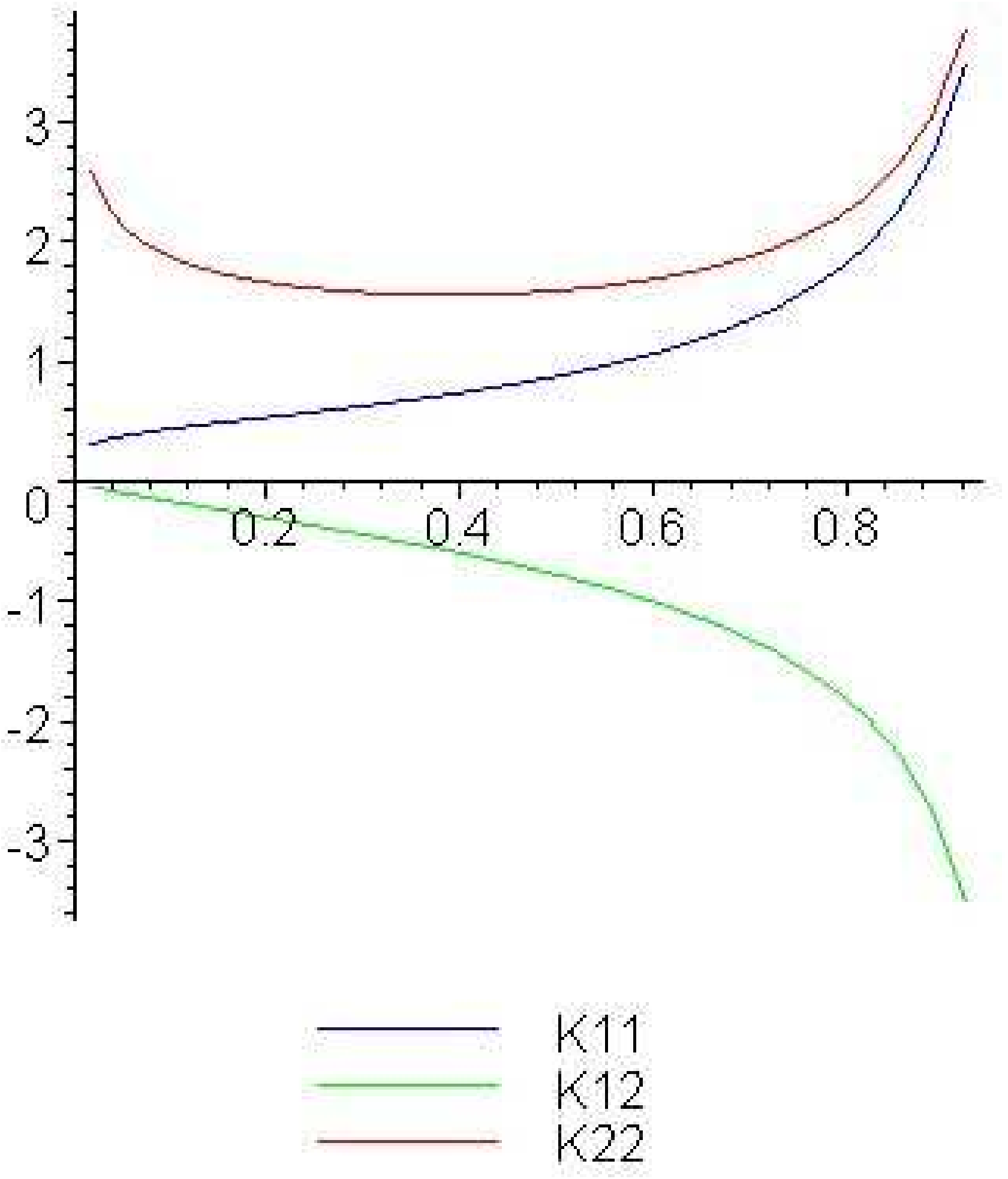}
\caption{$K$, $k_{0}=0.3$, $k_{1}=0.1$.}
\label{K31}
\end{figure}

\begin{figure}[t]
\centering
\includegraphics[width=2.6in]{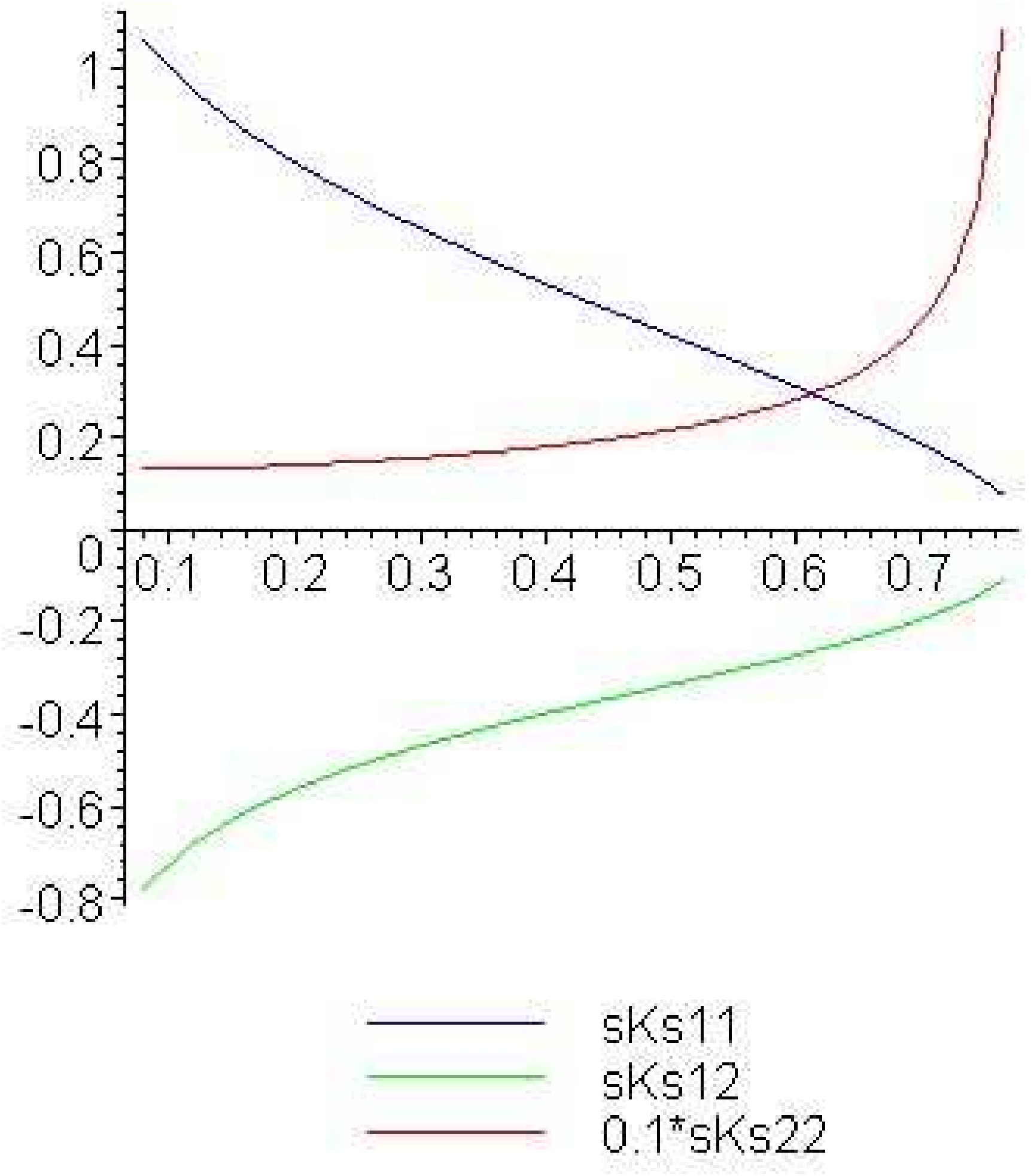}
\caption{$\sigma K\sigma$, $k_{0}=0.3$, $k_{1}=0.1$.}\label{sKs}

\end{figure}

We illustrate $K$ for $k_{0}=0.3$, $k_{1}=0.1$ with plots of $K (
\cos\theta,\sin\theta )  $ for $0<\theta<\frac{\pi}{4}$. Fig.~\ref{K31}
shows the values of $K_{11}$, $K_{12}$, $K_{22}$. For behavior near $x_{1}=x_{2}$
introduce
\[
\sigma=\frac{1}{\sqrt{2}}\left[
\begin{matrix}%
1 & 1\\
1 & -1
\end{matrix}
\right]  , \qquad \sigma\sigma^+_{12}\sigma=\left[
\begin{matrix}%
1 & 0\\
0 & -1
\end{matrix}
\right]  .
\]
Fig.~\ref{sKs} displays $\sigma K\sigma$ (thus $ (  \sigma
K\sigma )  _{12}=\frac{1}{2} (  K_{11}-K_{22} )  $; the $ (
2,2 )  $-entry is rescaled by $0.1$).

The degenerate cases $k_{0}=0$ and $k_{1}=0$ (when the group aspect reduces to
$\boldsymbol{Z}_{2}\times\boldsymbol{Z}_{2}$) provide a small check on the
calculations: for $k_{0}=0$ the weight is
\[
K(x)  =c\left[
\begin{matrix}%
 \vert x_{2}/x_{1} \vert ^{k_{1}} & 0\\
0 &  \vert x_{1}/x_{2} \vert ^{k_{1}}
\end{matrix}
\right]  ,
\]
and for $k_{1}=0$ and by use of the quadratic transformation $F\big(
a,a+\frac{1}{2};2a+1;1-u^{2}\big)  =\big(  \frac{1+u}{2}\big)  ^{-2a}$
(for $u$ near $1$)(see \cite[15.4.17]{Olver/Lozier/Boisvert/Clark:2010}) we obtain%
\[
K(x)  =c\sigma\left[
\begin{matrix}%
 \vert x_{1}-x_{2} \vert ^{k_{0}} \vert x_{1}+x_{2} \vert
^{-k_{0}} & 0\\
0 &  \vert x_{1}-x_{2} \vert ^{-k_{0}} \vert x_{1}+x_{2}%
 \vert ^{k_{0}}
\end{matrix}
\right]  \sigma.
\]

An orthogonal basis for $\mathcal{P}_{V}$ for the Gaussian inner product can
be given in terms of Laguerre polynomials and the harmonic polynomials from
the previous section. Recall from equations~(\ref{deltaF}) and~(\ref{xsqform})
(specializing $N=2$ and $\gamma (  \kappa;\tau )  =0$) that
\begin{gather*}
\begin{split}
& \Delta_{\kappa}^{k}\big(   \vert x \vert ^{2m}f(x,t)
\big)      =4^{k} (  -m )  _{k} (  -m-n )  _{k} \vert
x \vert ^{2m-2k}f(x,t)  ,\\
& \big\langle  \vert x \vert ^{2m}f, \vert x \vert
^{2m}g\big\rangle _{\tau}     =4^{m}m! (  n+1 )  _{m} \langle
f,g \rangle _{\tau},
\end{split}
\end{gather*}
for $f,g\in\mathcal{H}_{V,\kappa,n}$. These relations are transferred to the
Gaussian inner product
\[
 \langle f,g \rangle _{\tau}=\big\langle
e^{-\Delta_{\kappa}/2}f,e^{-\Delta_{\kappa}/2}g\big\rangle _{G}
\] by
computing
\begin{gather*}
e^{-\Delta_{\kappa}/2} \vert x \vert ^{2m}f(x,t)
=\sum_{j=0}^{m}\left(  -\frac{1}{2}\right)  ^{j}\frac{4^{j}}{j!} (
-m )  _{j} (  -m-n )  _{j} \vert x \vert ^{2m-2j}
f(x,t) \\
\hphantom{e^{-\Delta_{\kappa}/2} \vert x \vert ^{2m}f(x,t)  }{}
  = (  -1 )  ^{m}2^{m}m!L_{m}^{ (  n )  }\left(
\frac{ \vert x \vert ^{2}}{2}\right)  f(x,t)  ,
\end{gather*}
for $f\in\mathcal{H}_{V,\kappa,n}$ and $m\geq0$; where $L_{m}^{ (
n )  }$ denotes the Laguerre polynomial of degree~$m$ and index~$n$
(orthogonal for $s^{n}e^{-s}ds$ on $\mathbb{R}_{+}$). Denote $ \langle
f,f \rangle _{G}$ by $\nu_{G}(f)  $; recall $\nu (
f )  = \langle f,f \rangle _{\tau}$ for~$f\in\mathcal{P}_{V}$.

\begin{proposition}
The polynomials $L_{m}^{ (  n )  }\big(  \frac{ \vert
x \vert ^{2}}{2}\big)  p_{n,i}(x,t)  $ for $m,n\geq0$ and
$1\leq i\leq4$ $($except $i=1,2$ when $n=0)$ are mutually orthogonal in
$ \langle \cdot,\cdot \rangle _{G}$ and $\nu_{G} (  L_{m}^{ (
n )  } (   \vert x \vert ^{2}/2 )  p_{n,i} (
x,t )   )  =\frac{ (  n+1 )  _{m}}{m!}\nu (
p_{n,i} )  $.
\end{proposition}

The factor with $\nu (  p_{n,i} )  $ results from a simple calculation.
Note $\nu_{G}(f)  =\nu(f)  $ for any harmonic~$f$.
Because here $\gamma (  \kappa;\tau )  =0$ the harmonic decomposition
formula~(\ref{Hdecomp}) is valid for any parameter values. This is a notable
dif\/ference from the scalar case $\tau=1$ where $\gamma (  \kappa;1 )
=2k_{0}+2k_{1}$ and $2k_{0}+2k_{1}\neq-1,-2,\ldots$ is required for validity.

Using the same arguments as in the scalar case (see \cite[Section~5.7]{Dunkl/Xu:2001}) we def\/ine the \textit{Fourier transform} (for suitably
integrable functions~$f$). To adapt $E(x,y)  $ to vector notation
write $E(x,y)  =\sum\limits_{i,j=1}^{2}E_{ij}(x,y)
s_{i}t_{j}$ and set%
\begin{gather*}
\begin{split}
& \mathcal{F}f(y)  _{l}     :=\int_{\mathbb{R}^{2}}\sum_{i,j=1}
^{2}E_{li} (  x,-\mathrm{i}y )  K(x)  _{ij}f_{j} (
x )  dx_{1}dx_{2},\qquad l=1,2,\\
& \mathcal{F}f(y)      :=\mathcal{F}f(y)  _{1}
s_{1}+\mathcal{F}f(y)  _{2}s_{2}.
\end{split}
\end{gather*}

For $m,n\geq0$ and $1\leq i\leq4$ let $\phi_{m,n,i}(x)
=L_{m}^{ (  n )  } \big(   \vert x \vert ^{2} \big)
p_{n,i}(x)  e^{- \vert x \vert ^{2}/2}$.

\begin{proposition}
Suppose $m,n\geq0$ and $1\leq i\leq4$ then $\mathcal{F}\phi_{m,n,i} (
y )  = (  -\mathrm{i} )  ^{m+2n}\phi_{m,n,i}(y)
^{\ast}$. If $f(x)  =e^{- \vert x \vert ^{2}/2}g (
x )  $ for $g\in\mathcal{P}_{V}$ then $\mathcal{F} (  \mathcal{D}%
_{j}f )  (y)  =\mathrm{i}y_{j}\mathcal{F}f(y)
$, for $j=1,2$.
\end{proposition}

This establishes a Plancherel theorem for $\mathcal{F}$ by use of the density
(from Hamburger's theorem) of $\mathrm{span} \{  \phi_{m,n,i} \}  $
in $L^{2} (  K(x)  dx,\mathbb{R}^{2} )  $.

\section{Closing remarks}\label{section6}

The well-developed theory of the hypergeometric function allowed us to f\/ind
the weight function which satisf\/ies both a dif\/ferential equation and geometric
conditions. The analogous problem can be stated for any real ref\/lection group
and there are some known results about the dif\/ferential system~(\ref{eqnL})
(see \cite{Dunkl1993a,Dunkl1993b}); it appears some new insights are
needed to cope with the geometric conditions. The fact that the Gaussian inner
product $ \langle \cdot,\cdot \rangle _{G}$ is well-def\/ined supports
speculation that Gaussian-type weight functions exist in general settings.
However it has not been shown that~$K$ can be produced as a product~$L^{T}L$,
and the ef\/fect of the geometry of the mirrors (walls) on the solutions of the
dif\/ferential system is subtle, as seen in the $B_{2}$-case.

\subsection*{Acknowledgements}

This is the expanded version of an invited lecture presented at the Conference on Harmonic Analysis, Convolution Algebras, and Special Functions, TU
M\"{u}nchen, September~10,~2012.

\pdfbookmark[1]{References}{ref}
\LastPageEnding


\begin{thebibliography}{99}
\footnotesize\itemsep=0pt

\bibitem{Carter1993}
Carter R.W., Finite groups of {L}ie type. Conjugacy classes and complex
  characters, \textit{Wiley Classics Library}, John Wiley \& Sons Ltd., Chichester,
  1993.

\bibitem{Dunkl1993a}
Dunkl C.F., Dif\/ferential-dif\/ference operators and monodromy representations of
  {H}ecke algebras, \textit{Pacific~J. Math.} \textbf{159} (1993), 271--298.

\bibitem{Dunkl1993b}
Dunkl C.F., Monodromy of hypergeometric functions for dihedral groups,
  \href{http://dx.doi.org/10.1080/10652469308819011}{\textit{Integral Transform. Spec. Funct.}} \textbf{1} (1993), 75--86.

\bibitem{Dunkl/Opdam:2003}
Dunkl C.F., Opdam E.M., Dunkl operators for complex ref\/lection groups,
  \href{http://dx.doi.org/10.1112/S0024611502013825}{\textit{Proc. London Math. Soc.}} \textbf{86} (2003), 70--108,
  \href{http://arxiv.org/abs/math.RT/0108185}{math.RT/0108185}.

\bibitem{Dunkl/Xu:2001}
Dunkl C.F., Xu Y., Orthogonal polynomials of several variables,
  \href{http://dx.doi.org/10.1017/CBO9780511565717}{\textit{Encyclopedia of Mathematics and its Applications}}, Vol.~81, Cambridge
  University Press, Cambridge, 2001.

\bibitem{Etingof/Stoica:2009}
Etingof P., Stoica E., Unitary representations of rational {C}herednik
  algebras, \href{http://dx.doi.org/10.1090/S1088-4165-09-00356-2}{\textit{Represent. Theory}} \textbf{13} (2009), 349--370,
  \href{http://arxiv.org/abs/0901.4595}{arXiv:0901.4595}.

\bibitem{Griffeth2010}
Grif\/feth S., Orthogonal functions generalizing {J}ack polynomials,
  \href{http://dx.doi.org/10.1090/S0002-9947-2010-05156-6}{\textit{Trans. Amer. Math. Soc.}} \textbf{362} (2010), 6131--6157,
  \href{http://arxiv.org/abs/0707.0251}{arXiv:0707.0251}.

\bibitem{Olver/Lozier/Boisvert/Clark:2010}
Olver F.W.J., Lozier D.W., Boisvert R.F., Clark C.W. (Editors), N{IST} handbook
  of mathematical functions, U.S. Department of Commerce National Institute of
  Standards and Technology, Washington, DC, 2010.

\end{thebibliography}
\end{document}